\documentclass[a4paper]{article}
\usepackage{amsfonts,amscd}

\newcommand{\bea}{\begin{eqnarray*}}
\newcommand{\eea}{\end{eqnarray*}}
\newcommand{\bthm}[1]{\begin{thms}\label{t#1}}
\newcommand{\bd}[1]{\begin{defns}\label{d#1}}
\newcommand{\bprp}[1]{\begin{props}\label{p#1}}
\newcommand{\br}[1]{\begin{rems}\label{r#1}}
\newcommand{\bl}[1]{\begin{lems}\label{l#1}}
\newcommand{\ethm}{\end{thms}}
\newcommand{\ed}{\end{defns}}
\newcommand{\eprp}{\end{props}}
\newcommand{\er}{\end{rems}}
\newcommand{\el}{\end{lems}}
\newcommand{\lrf}[1]{Lemma~\ref{l#1}}
\newcommand{\trf}[1]{Theorem~\ref{t#1}}
\newcommand{\rrf}[1]{Remark~\ref{r#1}}
\newcommand{\prf}[1]{Proposition~\ref{p#1}}

\newcommand{\srf}[1]{Section~\ref{s#1}}

\newcommand{\bq}[1]{\begin{equation}\label{#1}}
\newcommand{\qq}{\end{equation}}
\newcommand{\bco}[1]{\begin{coros}\label{c#1}}
\newcommand{\crf}[1]{Corollary~\ref{c#1}}
\newcommand{\eco}{\end{coros}}
\newcommand{\een}{\end{enumerate}}

\newcommand{\beq}{\begin{equation}}
\newcommand{\eeq}{\end{equation}}
\newcommand{\rf}[1]{(\ref{#1})}
\newcommand{\beea}{\begin{eqnarray}}
\newcommand{\eeaa}{\end{eqnarray}}
\newcommand{\brr}[1]{\begin{equation}\label{#1}\begin{array}{ll}}
\newcommand{\err}{\end{array}\end{equation}}
\newcommand{\barr}{\begin{array}}
\newcommand{\earr}{\end{array}}
\def\ZZ{{\mathbb Z}}

\newtheorem{thms}{THEOREM}
\newtheorem{lems}{LEMMA}
\newtheorem{rems}{REMARK}
\newtheorem{defns}{DEFINITION}
\newtheorem{props}{PROPOSITION}
\newtheorem{coros}{COROLLARY}

\begin{document}

 \title{Complexity of Janet basis of a $D$-mo\-du\-le}
\author{Alexander Chistov
\\[-1pt]
%Dima Grigoriev
\small Steklov Institute of Mathematics,\\[-3pt]
\small Fontanka 27, St. Petersburg 191023, Russia\\[-3pt]
{\tt \small alch@pdmi.ras.ru}\\[-3pt]
%\\[-1pt]
%\small IRMAR, Universit\'e de Rennes \\[-3pt]
%\small Beaulieu, 35042, Rennes, France\\[-3pt]
%{\tt \small dima@math.univ-rennes1.fr}\\[-3pt]
%\small http://name.math.univ-rennes1.fr/dimitri.grigoriev
\and
Dima Grigoriev\\[-1pt]
\small CNRS, IRMAR, Universit\'e de Rennes \\[-3pt]
\small Beaulieu, 35042, Rennes, France\\[-3pt]
{\tt \small dmitry.grigoryev@univ-rennes1.fr}\\[-3pt]
\small http://perso.univ-rennes1.fr/dmitry.grigoryev
%\small Steklov Institute of Mathematics,\\[-3pt]
%\small Fontanka 27, St. Petersburg 191011, Russia\\[-3pt]
%{\tt \small inp@pdmi.ras.ru}\\[-3pt]
%\thanks{Partially supported by RFFI, grant 02-01-00093}
%\small http://www.pdmi.ras.ru/\~{}inp
}

\date{}

 \maketitle

\begin{abstract}
We prove a double-exponential upper bound on the degree and on the
complexity of
constructing a Janet basis of a $D$-module. This generalizes a well known
bound on
the complexity of a Gr\"obner basis of a module over the algebra of
polynomials. We would like to emphasize that the obtained bound
can not be immediately deduced from the commutative case.
\end{abstract}

\section*{Introduction}
Let $A$ be the Weyl algebra $F[X_1,\dots , X_n,$
${\partial \over \partial X_1},$ $ \dots ,
{\partial \over \partial X_n} ]$ (or
the algebra of differential operators $F(X_1,$ $\dots , X_n)
[{\partial \over \partial X_1}, \dots ,
{\partial \over \partial X_n} ]$). Denote for brevity
$D_i={\partial\over\partial X_i}$,
$1\le i\le n$. Any
$A$--module is called $D$--module.
It is well known that an $A$--module which is a
submodule of a free finitely generated $A$-mo\-du\-le has a
Janet basis.
Historically, it was first introduced in \cite{Janet}. In more recent
times of developing computer algebra Janet bases were studied in
\cite{Galligo},
\cite{Schwarz}, \cite{Li}. Janet bases generalize Gr\"obner bases
which were widely elaborated
in the algebra of polynomials (see e.~g.\cite{Cox}). For Gr\"obner bases a
double-exponential complexity bound was obtained in \cite{Mora},
\cite{Giusti}
relying on \cite{Bayer} and which was made more precise (with a
self--contained proof)
in \cite{Dube}.

Surprisingly, no complexity bound on Janet bases was established so far; in
the
present paper we fill this gap and prove a double-exponential complexity
bound.
On the other hand, a double-exponential complexity lower bound on Gr\"obner
bases
\cite{Mora}, \cite{Yap} provides by the same token a bound on Janet bases.

There is a folklore opinion that the problem of constructing a Janet basis
is easily reduced to the commutative case by considering the
associated graded mo\-du\-le, and, on the other hand, in the commutative
case \cite{Giusti}, \cite{Mora}, \cite{Dube} the
double--exponential upper bound is well known.
{\it But it turns out to be a fallacy! From a known
system of generators of a $D$-mo\-du\-le one can not obtain
immediately any system of
generators (even not necessarily a Gr\"obner basis) of
the associated
graded mo\-du\-le.} The main problem here is to construct such a
system of generators of the graded mo\-du\-le. It may have the
elements of degrees $(dl)^{2^{O(n)}}$, see the notation below.
Then, indeed, to the last system of generators of big degrees one can apply
the result known in the commutative case and get the bound
$((dl)^{2^{O(n)}})^{2^{O(n)}}=(dl)^{2^{O(n)}}$. So new ideas specific
to non--commutative case are needed.

We are interested in the estimations for Janet bases of $A$-submodules of
$A^l$. The Janet basis depends on the choice of the linear order on the
monomials (we define them also for $l>1$).
In this paper we consider the most
general linear orders on the monomials from
$A^l$. They satisfy conditions (a) and (b) from \srf{1} and are called {\it
admissible}.
We prove the following result.

\bthm{1} For any admissible linear order on the monomials from
$A^l$ any $A$-submodule $I$ of $A^l$
generated by elements of degrees at
most $d$ (with respect to the filtration in the corresponding algebra, see
Section~\ref{s1} and \srf{8})
has a Janet basis with the degrees and the number of
its elements less than
\[
(dl)^{2^{O(n)}}.
\]
\end{thms}
%metka t1 \\

\noindent We prove in detail
this theorem for the case of the Weyl al\-ge\-b\-ra $A$.
The proof for the case of the al\-ge\-b\-ra of differential operators is
similar.
It is sketched in \srf{8}.
>From \trf{1} we get that  the Hilbert
function $H(I,m)$, see \srf{1}, of the $A$-submodule from
this theorem is stable for $m\ge(dl)^{2^{O(n)}}$ and
the absolute values of all coefficients of the Hilbert
polynomial of $I$ are bounded from above by $(dl)^{2^{O(n)}}$,  cf. e.g.,
\cite{Mora}.
This fact follows directly from  \rf{65}, \lrf{14} from Appendix~1,
\lrf{3} and \trf{2}.
We mention that in \cite{G05} the similar bound was
shown on the leading coefficient of the Hilbert polynomial.

\medskip Now we outline the plan for the proof of \trf{1}.
The main tool in the proof is a homogenized Weyl algebra ${^h\!A}$ (or
respectively, a homogenized algebra of differential operators ${^h\!B}$).
It is introduced in
Section~\ref{s3} (respectively, Section~\ref{s8}).
The al\-ge\-b\-ra ${^h\!A}$ (respectively ${^h\!B}$) is generated over the
ground
field $F$ by
$X_0,\ldots ,X_n$, $D_1,\ldots , D_n$
(respectively over the field $F(X_1,\ldots ,X_n)$
by $X_0,D_1,\ldots , D_n$). Here $X_0$ is a new homogenizing variable.
In the al\-ge\-b\-ra ${^h\!A}$ (respectively ${^h\!B}$)
relations \rf{8} \srf{3} (respectively \rf{63} \srf{8}) hold for
these generators in ${^h\!A}$.

We define the homogenization ${^h\!I}$ of the mo\-du\-le $I$. It is a
${^h\!A}$--submodule of
${^h\!A^l}$. The main problem is to estimate the degrees of a system of
generators of
${^h\!I}$. These estimations are central in the paper. They are
deduced from \trf{2} \srf{5}. This theorem is devoted to the problem
of solving systems of linear equations over the ring ${^h\!A}$; we discuss
it below in more detail.

The system of generators of ${^h\!I}$ gives a system of generators of the
graded
$\mathop{\rm gr}\nolimits(A)$--mo\-du\-le $\mathop{\rm gr}\nolimits(I)$
corresponding to $I$. But $\mathop{\rm gr}\nolimits(A)$ is a polynomial ring.
Hence using \lrf{14} Appendix~1 we get a double--exponential
bound $(dl)^{2^{O(n)}}$ on the stabilization of the Hilbert
function of $\mathop{\rm gr}\nolimits(I)$ and
the absolute values of the coefficients of the
Hilbert polynomial of $\mathop{\rm gr}\nolimits(I)$.
Therefore, the similar bound holds for the stabilization of the Hilbert
functions of $I$ and the coefficients of the Hilbert polynomial of $I$, see
\srf{2}.

But the Hilbert functions of the mo\-du\-les $I$ and ${^h\!I}$ coincide, see
\srf{3}.
Hence the last bound holds also for the stabilization of the Hilbert
functions of  ${^h\!I}$ and the coefficients of the Hilbert polynomial of
${^h\!I}$.
In \srf{4} we introduce the linear order on the monomials from ${^h\!A^l}$
induced by the initial linear order on the monomials from $A^l$
(the homogenizing variable $X_0$ is the least possible in this ordering).
Further, we define the Janet basis of ${^h\!I}$
with respect to the induced linear order on the monomials. Such a basis can
be obtained by the homogenization
of the elements of a Janet basis of $I$ with respect to the initial linear
order, see \lrf{15}.

Let $\mathop{\rm
Hdt}\nolimits({^h\!I})$ be the monomial
mo\-du\-le (i.e., the mo\-du\-le
which has a system of generators consisting of monomials)
generated by the greatest monomials of all the elements of
the mo\-du\-le ${^h\!I}$, see \srf{10}.
Let ${^c\!I}$, see \srf{10}, be the mo\-du\-le over the polynomial ring
${^c\!A}=F[X_0,\ldots , X_n,D_1,\ldots , D_n]$
generated by all the monomials
from $\mathop{\rm Hdt}\nolimits({^h\!I})$ (they are
considered now as elements of ${^c\!A}$).
Then the Hilbert functions of the mo\-du\-les
${^h\!I}$ and ${^c\!I}$ coincide.
Thus, we have the same as above double--exponential estimation for
the stabilization of the Hilbert
functions of ${^c\!I}$ and the coefficients of the
Hilbert polynomial of ${^c\!I}$. Now using \lrf{16}
we get the estimation $(dl)^{2^{O(n)}}$
on the monomial system of generators of ${^c\!I}$,
hence also of
$\mathop{\rm Hdt}\nolimits({^h\!I})$. This gives the bound for the degrees
of the elements of the
Janet bases of ${^h\!I}$ and hence also for the required Janet basis of $I$,
and proves \trf{1}.

\medskip
The problem of solving systems of linear equations over the
homogenized algebra is central in this paper, see \trf{2}.
It is studied in Sections~\ref{s4}--\ref{s5}. A similar problem over
the Weyl algebra (without a homogenization) was considered in \cite{G05}. The
principal idea is to try to extend the well known method due to G.Hermann
\cite{Herm} which was elaborated for the algebra of
polynomials, to
the homogenized Weyl algebra.
There are two principal difficulties on this way.
The first one is that in the method of
G.Hermann the use of determinants is essential which one has to avoid
dealing with non-commutative algebras. The second is that
one needs a kind of the Noether normalization theorem
in the situation under consideration. So it is necessary
to choose the leading elements in the analog of the G.Hermann method
with the least $\mathop{\rm ord}\nolimits_{X_0}$, where $X_0$ is a
homogenizing variable, see \srf{3}.

The obtained bound on the degree of a Janet basis
implies a similar bound on the complexity of its constructing.
Indeed, by \crf{1} (it is formulated for the case of Weyl al\-ge\-b\-ra but
the analogous corollary holds for the case of al\-ge\-b\-ra of differential
operators) one can compute the linear space of all the elements
$z\in I$ of degrees bounded from above by $(dl)^{2^{O(n)}}$. Further,
by \trf{1} the mo\-du\-le $\mathop{\rm Hdt}\nolimits(I)$,
see \srf{1}, is generated by all
the elements $\mathop{\rm Hdt}\nolimits(z)$ with $z\in I$ of
degrees bounded from above by $(dl)^{2^{O(n)}}$.
Hence one can compute a system of
generators of $\mathop{\rm Hdt}\nolimits(I)$ and a Janet basis of $I$
solving linear systems over $F$ of size bounded from above $(dl)^{2^{O(n)}}$
(just by the enumeration of all monomials
of degrees at most $(dl)^{2^{O(n)}}$ which are possible generators of
$\mathop{\rm Hdt}\nolimits(I)$).
If one needs to construct the reduced Janet basis it is sufficient to apply
additionally \rrf{3} \srf{10}.

For the sake of self--containedness in Appendix~1, see \lrf{14},
we give a short proof of
the double--exponential estimation for stabilization of
the Hilbert function of a graded mo\-du\-le over a homogeneous polynomial
ring.
%It is well--known that from this one can deduce the estimation for degrees
%of elements of a Gr\"obner basis of an ideal (if the considered mo\-du\-le
%is an ideal), see, e.g., \cite{Giusti},~\cite{Mora}.
A conversion of \lrf{14} also holds, see Appendix~1 \lrf{16}.
It is essential for us. The proof of \lrf{16} uses the classic
description of the Hilbert function of a homogeneous ideal
in $F[X_0,\ldots , X_n]$ via Macaulay
constants $b_{n+2},\ldots , b_1$ and the constant $b_0$ introduced in
\cite{Dube}. In Appendix~2 we give an independent and instructive proof
of \prf{monomial} which is similar to \lrf{16}.
In some sence \prf{monomial} is even more
strong than \lrf{16} since to apply it one does not need a bound
for the stabilization of the Hilbert function. Of course, the reference to
\prf{monomial} can be used in place of \lrf{16} in our paper.
%To give a reference
%we tried to find a result similar to \lrf{16}
%in literature but did not succeed.

\section{Definition of the Janet basis}\label{s1}
%metka s1 \\

\noindent Let $A=F[X_1,\ldots,X_n,D_1,\ldots , D_n]$, $n\ge 1$, be a Weyl
al\-ge\-b\-ra over a field $F$ of
zero--characteristic. So $A$ is defined by the following relations
\bq{7}
X_vX_w=X_wX_v,\;D_vD_w=D_wD_v,\;D_vX_v-X_vD_v=1,\;X_vD_w=D_wX_v,\quad v\ne w.
\end{equation}
%metka 7 \\
By  \rf{7} any element $f\in A$ can be uniquely represented in the form
\bq{10}
f=\sum_{i_1,\ldots , i_n,j_1,\ldots , j_n\ge 0}f_{i_1,\ldots ,
i_n,j_1,\ldots , j_n}
X_1^{i_1}\ldots X_n^{i_n}D_1^{j_1}\ldots D_n^{j_n},
\end{equation}
%metka 10 \\
where all $f_{i_1,\ldots , i_n,j_1,\ldots , j_n}\in F$ and only a finite
number of $f_{i_1,\ldots , i_n,j_1,\ldots , j_n}$ are nonzero. Denote for
brevity  ${\mathbb Z}_+=\{z\in{\mathbb Z}\, :\, z\ge 0\}$ to the set of all
nonnegative integers and
\bq{11}
\begin{array}{l}
i=(i_1,\ldots , i_n),\quad j=(j_1,\ldots , j_n),\quad
f_{i,j}=f_{i_1,\ldots , i_n,j_1,\ldots , j_n} \\
X^i=X_1^{i_1}\ldots X_n^{i_n},\quad D^j=D_1^{j_1}\ldots D_n^{j_n},\quad
f=\sum_{i,j}f_{i,j}X^iD^j,\\
|i|=i_1+\ldots + i_n,\quad i+j=(i_1+j_1,\ldots , i_n+j_n).
\end{array}
\end{equation}
%metka 11 \\
So $i,j\in{\mathbb Z}_+^n$ are multiindices.
By definition the degree of $f$
\[
\deg f=\deg_{X_1,\ldots , X_n,D_1,\ldots , D_n}f=\max\{|i|+|j|\, :\,
f_{i,j}\ne 0\}.
\]
Let $M$ be a left $A$-mo\-du\-le given by its
generators
$m_1,\ldots, m_l$, $l\ge 0$, and relations
\bq{1}
\sum_{1\le w\le l}a_{v,w}m_w,\quad 1\le v\le k.
\end{equation}
%metka 1 \\
where $k\ge 0$ and
all $a_{v,w}\in A$.
We assume that $\deg a_{v,w}\le d$ for all $v,w$.
By \rf{1} we have the exact sequence
\bq{3}
A^k\stackrel{\iota}{\rightarrow} A^l\stackrel{\pi}{\rightarrow}
M\rightarrow 0
\end{equation}
%metka 3 \\
of left $A$-mo\-du\-les. Denote
$I=\iota(A^k)\subset A^l$. If $l=1$ then $I$ is
a left
ideal of $A$ and $M=A/I$.
In the general case $I$ is generated by the elements
\[
(a_{v,1},\ldots ,
a_{v,l})\in A^l,\quad 1\le v\le k.
\]
For an integer $m\ge 0$ put
\bq{4}
A_m=\{a\, :\,\deg a\le m\},\quad M_m=\pi(A_m^l),\quad I_m=I\cap A_m^l.
\end{equation}
%metka 4 \\
So now $A$, $M$, $I$ are filtered mo\-du\-les with filtrations $A_m$, $M_m$,
$I_m$,
$m\ge 0$, respectively and the sequence of ho\-mo\-mor\-ph\-isms of
vector spaces
\[
0\rightarrow I_m\rightarrow A_m^l\rightarrow M_m\rightarrow 0
\]
induced by  \rf{3} is exact for every $m\ge 0$. The Hilbert function
$H(M,m)$ of the mo\-du\-le
$M$
is defined by the equality
\[
H(M,m)=\dim_F M_m,\quad m\ge 0.
\]

Each element of $A^l$ can be uniquely represented as an $F$-linear
combination of elements
$e_{v,i,j}=(0,\ldots , 0, X^iD^j,0,\ldots , 0)$,
herewith $i,j\in{\mathbb Z}_+^n$ are multiindices, see \rf{11},
and the nonzero monomial  $X^iD^j$ is at the
position $v$, $1\le v\le
l$. So every element $f\in A^l$ can be represented in the form
\bq{2}
f=\sum_{v,i,j}f_{v,i,j}e_{v,i,j},\quad f_{v,i,j}\in F.
\end{equation}
%metka 2 \\
The elements $e_{v,i,j}$ will be called monomials.

Consider a linear order $<$ on the set of all the monomials
$e_{v,i,j}$ or which is the same
on the set of triples $(v,i,j)$, $1\le v\le l$, $i,j\in{\mathbb Z}_+^n$.
If $f\ne 0$ put
\bq{43}
o(f)=\max\{(v,i,j)\, :\, f_{v,i,j}\ne 0\},
\end{equation}
%metka 43 \\
see \rf{2}. Set
\[
o(0)=-\infty<o(f)
\]
for every $0\ne f\in A$.
Let us define the leading monomial of the element $0\ne f\in
A^l$ by the formula
\[
\mathop{\rm Hdt}\nolimits(f)=f_{v,i,j}e_{v,i,j},
\]
where $o(f)=(v,i,j)$. Put $\mathop{\rm
Hdt}\nolimits(0)=0$.
Hence $o(f-\mathop{\rm
Hdt}\nolimits(f))<o(f)$ if $f\ne 0$.
For $f_1,f_2\in A^l$ if $o(f_1)<o(f_2)$ we shall write $f_1<f_2$.
We shall require additionally that
\begin{enumerate} \renewcommand{\labelenumi}{(\alph{enumi})}
\item for all multiindices $i,j,i',j'$ for all $1\le v\le l$
if $i_1\le i'_1,\ldots , i_n\le i'_n$ and $j_1\le j'_1,\ldots , j_n\le j'_n$
then $(v,i,j)\le(v,i',j')$.
\item for all multiindices $i,j,i',j',i'',j''$ for all $1\le v,v'\le l$ \
if  $(v,i,j)<(v',i',j')$ then $(v,i+i'',j+j'')<(v',i'+i'',j'+j'')$.
\end{enumerate}
Conditions (a) and (b) imply that for all $f_1,f_2\in A^l$ for every
nonzero  $a\in A$
if $f_1<f_2$ then $af_1<af_2$, i.e., the considered linear order is
compatible with the
products. Any linear order on monomials $e_{v,i,j}$
satisfying (a) and (b) will be
called {\it admissible}.

Set
\[
\mathop{\rm Hdt}\nolimits(I)=
\sum_{f\in I}A\mathop{\rm Hdt}\nolimits(f).
\]
So $\mathop{\rm Hdt}\nolimits(I)$ is an ideal of $A$.
By definition the family $f_1,\ldots ,f_m$ of elements of $I$ is a Janet
basis
of the mo\-du\-le $I$ if and only if
\begin{enumerate} \renewcommand{\labelenumi}{\arabic{enumi})}
\item $\mathop{\rm Hdt}\nolimits(I)=A\mathop{\rm Hdt}\nolimits(f_1)+\ldots +
A\mathop{\rm Hdt}\nolimits(f_m)$, i.e., the submodule of $A^l$ generated by
$\mathop{\rm Hdt}\nolimits(f_1),\ldots,\mathop{\rm
Hdt}\nolimits(f_m)$ coincides with $\mathop{\rm
Hdt}\nolimits(I)$.
\end{enumerate}
Further, the Janet basis $f_1,\ldots ,f_m$ of $I$ is reduced if
and only if the following conditions hold.
\begin{enumerate} \renewcommand{\labelenumi}{\arabic{enumi})}
\setcounter{enumi}{1}
\item $f_1,\ldots ,f_m$ does not
contain a smaller Janet basis of $I$,
\item $\mathop{\rm Hdt}\nolimits(f_1)>\ldots >\mathop{\rm
Hdt}\nolimits(f_m)$.
\item the coefficient from $F$ of every monomial
$\mathop{\rm Hdt}\nolimits(f_v)$,
$1\le v\le l$, is $1$.
\item Let $f_\alpha=\sum_{v,i,j}f_{\alpha,v,i,j}e_{v,i,j}$ be representation
\rf{10} for $f_\alpha$, $1\le\alpha\le m$. Then for all
$1\le\alpha<\beta\le m$ for all $1\le v\le l$ and
multiindices $i,j$ the monomial
$f_{\alpha,v,i,j}e_{v,i,j}\not\in\mathop{\rm
Hdt}\nolimits(Af_\beta\setminus\{0\})$.
\end{enumerate}
Since the ring $A$ is Noetherian for considered $I$ there exists a
Janet basis. Further the reduced Janet basis of $I$ is
uniquely defined.

\section{The graded
mo\-du\-le corresponding to a $D$--mo\-du\-le }\label{s2}
%metka s2 \\

Put $A_v=I_v=M_v=0$ for $v<0$ and
\[
\mathop{\rm gr}\nolimits (A)=\oplus_{m\ge 0}A_m/A_{m-1},\;\mathop{\rm
gr}\nolimits(I)=\oplus_{m\ge 0}I_m/I_{m-1},\;
\mathop{\rm gr}\nolimits(M)=\oplus_{m\ge 0}M_m/M_{m-1}.
\]
The structure of the al\-ge\-b\-ra on $A$ induces the structure  of a graded
al\-ge\-b\-ra on $\mathop{\rm gr}\nolimits(A)$. So we have
$\mathop{\rm gr}\nolimits(A)=F[X_1,\ldots ,X_n,D_1,\ldots , D_n ]$ is an
al\-ge\-b\-ra of
polynomials with respect to the variables
$X_1,\ldots ,X_n$, $D_1,\ldots ,D_n$. Further, $\mathop{\rm gr}\nolimits(I)$
and
$\mathop{\rm gr}\nolimits(M)$ are graded $\mathop{\rm
gr}\nolimits(A)$-mo\-du\-les. From  \rf{4} we get the exact sequences
\bq{5}
0\rightarrow I_m/I_{m-1}\rightarrow (A_m/A_{m-1})^l\rightarrow
M_m/M_{m-1}\rightarrow 0,\quad m\ge 0.
\end{equation}
%metka 5 \\
The Hilbert function of the mo\-du\-le $\mathop{\rm gr}\nolimits(M)$ is
defined as follows
\[
H(\mathop{\rm gr}\nolimits(M),m)=\dim_FM_m/M_{m-1},\quad m\ge 0.
\]
Obviously
\bq{65}
H(M,m)=\sum_{0\le v\le m}H(\mathop{\rm gr}\nolimits(M),v),\quad
H(\mathop{\rm gr}\nolimits(M),m)=H(M,m)-H(M,m-1).
\end{equation}
%metka 65 \\
for every $m\ge 0$.

Denote for an arbitrary $a\in M$ by $\mathop{\rm
gr}\nolimits(a)\in\mathop{\rm gr}\nolimits(M)$ the image of $a$ in
$\mathop{\rm gr}\nolimits(M)$.

\bl{2} Assume that $b_1,\ldots , b_s$ is a system of generators of $I$.
Let $\nu_i=\deg b_i$, $1\le i\le s$. Suppose that for every $m\ge 0$
\bq{6}
I_m=\Bigl\{\sum_{1\le v\le \mu}c_vb_v\, :\,c_v\in A,\quad\deg c_v\le
m-\nu_v,\quad
1\le i\le s \Bigr\}.
\end{equation}
%metka 6 \\
Then $\mathop{\rm gr}\nolimits(b_1),\ldots , \mathop{\rm
gr}\nolimits(b_s)$ is a system of generators of the $\mathop{\rm
gr}\nolimits(A)$-mo\-du\-le $\mathop{\rm gr}\nolimits(I)$.
\end{lems}
%metka l2 \\
\noindent{\bf PROOF}\quad This is straightforward.

\medskip\noindent So it is sufficient to construct a system of generators
$b_1,\ldots , b_s$ of $I$ satisfying \rf{6}.

\section{Homogenization of the Weyl al\-ge\-b\-ra}\label{s3}
%metka s3 \\

Let $X_0$ be a new variable. Consider the al\-ge\-b\-ra
${^h\!A}=F[X_0,X_1,\ldots
,X_n,D_1,$ $\ldots , D_n]$ given by the relations
\bq{8}
\begin{array}{l}
X_vX_w=X_wX_v,\;D_vD_w=D_wD_v,\quad\mbox{for all}\quad
 v,w, \\
D_vX_v-X_vD_v=X_0^2,\; 1\le v\le n,\quad X_vD_w=D_wX_v\quad\mbox{for
all}\quad
 v\ne w.\\
\end{array}
\end{equation}
%metka 8 \\
The al\-ge\-b\-ra ${^h\!A}$ is
Noetherian similarly to the Weyl al\-ge\-b\-ra $A$.
By  \rf{8} an element $f\in {^h\!A}$ can be uniquely represented in the form
\bq{15}
f=\sum_{i_0,i_1,\ldots , i_n,j_1,\ldots , j_n\ge 0}f_{i_0,\ldots ,
i_n,j_1,\ldots , j_n}
X_0^{i_0}\ldots X_n^{i_n}D_1^{j_1}\ldots D_n^{j_n},
\end{equation}
%metka 15 \\
where all $f_{i_0,\ldots , i_n,j_1,\ldots , j_n}\in F$ and only a finite
number of $f_{i_0,\ldots , i_n,j_1,\ldots , j_n}$ are nonzero.
Let $i,j$ be multiindices, see \rf{11}.
Denote for
brevity
\bq{12}
\begin{array}{l}
i=(i_1,\ldots , i_n),\quad j=(j_1,\ldots , j_n),\quad
f_{i_0,i,j}=f_{i_0,\ldots , i_n,j_1,\ldots , j_n} \\
f=\sum_{i_0,i,j}f_{i_0,i,j}X_0^{i_0}X^iD^j.
\end{array}
\end{equation}
%metka 12 \\
By definition the degrees of $f$
\begin{eqnarray*}
&&\deg f=\deg_{X_0,\ldots , X_n,D_1,\ldots , D_n}f=\max\{i_0+|i|+|j|\, :\,
f_{i_0,i,j}\ne 0\}, \\
&&\deg_{D_1,\ldots , D_n}f=\max\{|j|\, :\,
f_{i_0,i,j}\ne 0\}, \\
&&\deg_{D_\alpha}f=\max\{j_\alpha\, :\,
f_{i_0,i,j}\ne 0\},\quad 1\le \alpha\le n\\
&&\deg_{X_\alpha}f=\max\{i_\alpha\, :\,
f_{i_0,i,j}\ne 0\},\quad 1\le \alpha\le n
\end{eqnarray*}
Set $\mathop{\rm ord}\nolimits 0=\mathop{\rm ord}\nolimits_{X_0}0=+\infty$.
If $0\ne f\in{^h\!A}$
then put
\bq{44}
\mathop{\rm ord}\nolimits f=\mathop{\rm ord}\nolimits_{X_0} f
=\mu\quad\mbox{if and only if}\quad f\in
X_0^{\mu}({^h\!A})\setminus X_0^{\mu+1}({^h\!A}),\quad\mu\ge 0.
\end{equation}
%metka 44 \\
For every $z=(z_1,\ldots , z_l )\in{^h\!A^l}$ put
\[
\mathop{\rm
ord}\nolimits z=\min_{1\le i\le l}\{\mathop{\rm
ord}\nolimits z_i\}, \quad  \deg z=\max_{1\le i\le l}\{\deg z_i\}.
\]
Similarly one defines $\mathop{\rm ord}\nolimits b$ and $\deg b$ for an
arbitrary $(k\times l)$--ma\-t\-rix  $b$
with coefficients from ${^h\!A}$. More precisely, one consider here $b$ as a
vector with $kl$ entries.

The element $f\in{^h\!A}$ is homogeneous if and only if
$f_{i_0,i,j}\ne 0$ implies $i_0+|i|+|j|=\deg f$, i.e., if and only if
$f$ is a sum of monomials of the same
degree $\deg f$.
The homogeneous degree of a nonzero homogeneous element $f$ is its degree.
The homogeneous degree of $0$ is not defined ($0$ belongs to all the
homogeneous components of ${^h\!A}$, see below).

The $m$-th homogeneous component of ${^h\!A}$ is the $F$-linear space
\[
({^h\!A})_m=\left\{\,z\in {^h\!A}\,:\,z\;\mbox{is homogeneous}\;\&\;\deg
z=m\;\mbox{or}\; z=0\,\right\}
\]
for every integer $m$.
Now ${^h\!A}$ is a graded ring with respect to the homogeneous
degree.
By definition the ring ${^h\!A}$ is a homogenization of the Weyl
al\-ge\-b\-ra $A$.

We shall consider the category of finitely generated graded
mo\-du\-les $G$ over the ring ${^h\!A}$. Such a mo\-du\-le
$G=\oplus_{m\ge m_0}G_m$ is
a direct sum of its homogeneous components $G_m$, where $m,m_0$. are
integers. Every $G_m$ is a finite dimensional $F$-linear space and
$({^h\!A})_pG_m\subset G_{p+m}$
for all integers $p,m$. If $G$ and $G'$ are two finitely generated graded
${^h\!A}$-mo\-du\-les
then $\varphi\,:\, G\rightarrow G'$ is a morphism (of degree $0$) of the
graded
mo\-du\-les
if and only if $\varphi$ is a morphism of ${^h\!A}$-mo\-du\-les and
$\varphi(G_m)\subset G'_m$ for every integer $m$.

The element $z\in {^h\!A}$ (respectively $z\in A$)
is called to be the term if and only if $z=\lambda
z_1\cdot\ldots\cdot
z_\nu$ for some $0\ne\lambda\in F$, integer $\nu\ge 0$ and
$z_w\in\{X_0,\ldots ,X_n,D_1,\ldots , D_n\}$ (respectively
$z_w\in\{X_1,\ldots ,X_n,D_1,\ldots , D_n\}$),
$1\le w\le\nu$.

Let $z=\sum_{j}z_j\in A$ be an arbitrary element of  the Weyl
al\-ge\-b\-ra $A$ represented as a sum
of terms $z_j$ and $\deg z=\max_{j}\deg z_j$. One can take here,
for example, representation \rf{11} for $z$.
Then we define the homogenization ${^h\!z}\in{^h\!A}$ by the formula
\[
{^h\!z}=\sum_{j}z_jX_0^{\deg z-\deg z_j}.
\]
By \rf{7}, \rf{8} the right part of the last equality does not depend on the
chosen representation of
$z$ as a sum of terms. Hence ${^h\!z}$ is defined correctly.
If $z\in{^h\!A}$ then  ${^a\!z}\in A$ is obtained by substituting  $X_0=1$
in $z$.
Hence for every $z\in A$ we have ${^a\!{^h\!z}}=z$, and for every
$z\in{^h\!A}$ the
element $z={^h\!{^a\!z}}X_0^\mu$, where $\mu=\mathop{\rm
ord}\nolimits z$.

For an element $z=(z_1,\ldots , z_l)\in A^l$ put
$\deg z=\max_{1\le i\le l}\{\deg z_i\}$ and
\[
{^h\!z}=\left(\,{^h\!z}_1X_0^{\deg z-\deg z_1},\ldots
,{^h\!z}_lX_0^{\deg z-\deg z_l}\,\right)\in{^h\!A^l}.
\]
Similarly one defines $\deg a$ and the homogenization
${^h\!a}=(a_{v,w})_{1\le v\le k,\,1\le w\le l}$ for an
arbitrary $k\times l$--ma\-t\-rix $a$
with coefficients from $A$. More precisely, one consider here $a$ as a
vector with $kl$ entries.
Hence if $b=(b_{v,w})_{1\le v\le k,\,1\le w\le l}={^h\!a}$
then $b_{v,w}={^h\!a}_{v,w}X_0^{\deg a-\deg a_{v,w}}$ for all $v,w$.

The $m$-th homogeneous component of ${^h\!A^l}$ is
\[
({^h\!A^l})_m=\left\{\,{^h\!z}\,:\,z\in A^l\;\&\;\deg z=m\;\mbox{or}\;
z=0\,\right\}
\]
For an $F$-linear subspace $X\subset A^l$ put ${^h\!X}$ to be
the least linear subspace of ${^h\!A^l}$ containing the set
$\{{^h\!z}\, :\, z\in
X\}$. If $X$ is a (finitely generated) $A$-submodule of $A^l$
then ${^h\!X}$ is a (finitely generated)
graded submodule of ${^h\!A^l}$. The graduation on ${^h\!X}$ is induced by
the one of ${^h\!A^l}$.

For an element $z=(z_1,\ldots , z_l)\in
{^h\!A^l}$ put ${^a\!z}=({^a\!z}_1,\ldots
,{^a\!z}_l)\in A^l$.
For a subset $X\subset{^h\!A^l}$ put ${^a\!X}=\{{^a\!z}\, :\, z\in
X\}\subset A^l$. If $X$ is a $F$-linear space then
${^a\!X}$ is also a $F$-linear space. If $X$ is a finitely generated graded
submodule of
${^h\!A^l}$ then ${^a\!X}$ is finitely generated  submodule of $A^l$.

Now ${^h\!I}$ is a graded submodule of ${^h\!A^l}$.
Further, ${^a\!{^h\!I}}=I$.
Let $({^h\!I})_m$ be the $m$-th homogeneous component of ${^h\!I}$.
Then
\begin{eqnarray}
&&{^h\!(I_m)} =\oplus_{0\le j\le m}({^h\!I} )_j,\quad
m\ge 0, \label{28} \\
&&{^a\!(({^h\!I} )_m)}=I_m,\quad
m\ge 0. \label{29}
\end{eqnarray}
%metka 28,29, \\
and  \rf{29} induces the isomorphism $\iota\, :\,({^h\!I})_m\rightarrow
I_m$.
Set
${^h\!M}={^h\!A^l}/{^h\!I}$. Hence ${^h\!M}$ is a graded
${^h\!A}$-mo\-du\-le and we have the exact sequence
\bq{51}
0\rightarrow
{^h\!I}\rightarrow {^h\!A^l}\rightarrow {^h\!M}\rightarrow 0.
\end{equation}
%metka 51 \\
The $m$-th homogeneous component $({^h\!M})_m$ of ${^h\!M}$
\bq{30}
({^h\!M})_m=({^h\!A^l})_m/({^h\!I})_m\simeq A^l_m/I_m.
\end{equation}
%metka 30 \\
by the isomorphism $\iota$.
We have the exact sequences
\bq{52}
0\rightarrow ({^h\!I})_m\rightarrow({^h\!A^l})_m\rightarrow({^h\!M})_m
\rightarrow 0,\quad m\ge 0.
\end{equation}
%metka 52 \\
By definition the Hilbert function of the mo\-du\-le ${^h\!M}$ is
\[
H({^h\!M},m)=\dim_F({^h\!M})_m,\quad m\ge 0.
\]
By \rf{30} we have $H(M,m)=H({^h\!M},m)$ for every $m\ge 0$,
i.e., the Hilbert
functions of  $M$ and ${^h\!M}$ coincide.

\bl{3} Let $b_1,\ldots , b_s$ be a system of homogeneous generators of the
${^h\!A}$-mo\-du\-le ${^h\!I}$.
Then
\[
\mathop{\rm gr}\nolimits({^a\!b_1} ),\ldots ,\mathop{\rm
gr}\nolimits({^a\!b_s} )\in\mathop{\rm
gr}\nolimits(A)^l
\]
is a system of generators of $\mathop{\rm gr}\nolimits(A)$-mo\-du\-le
$\mathop{\rm gr}\nolimits(I)$.
\end{lems}
%metka l3 \\
\noindent{\bf PROOF}\quad By  \rf{29} $^a\!(({^h\!I)_m})=I_m$.
Now the required assertion follows from \lrf{2}. The lemma is proved.

\section{The Janet bases of a mo\-du\-le and of
its homogenization}\label{s10}
%metka s10 \\

Each element of ${^h\!A^l}$ can be uniquely represented as an $F$-linear
combination of elements
$e_{v,i_0,i,j}=(0,\ldots , 0, X_0^{i_0}X^iD^j,0,\ldots , 0)$,
herewith $0\le i_0\in{\mathbb Z}$, $i,j\in{\mathbb Z}_+^n$ are
multiindices, see \rf{11},
and the nonzero monomial  $X_0^{i_0}X^iD^j$ is at the
position $v$, $1\le v\le
l$. So every element $f\in {^h\!A^l}$ can be represented in the form
\bq{2a}
f=\sum_{v,i_0,i,j}f_{v,i_0,i,j}e_{v,i_0,i,j},\quad f_{v,i_0,i,j}\in F.
\end{equation}
%metka 2a \\
and only a finite number of $f_{v,i_0,i,j}$ are nonzero.
The elements $e_{v,i_0,i,j}$ will be called monomials.

Let us replace everywhere in \srf{1} after the definition of the Hilbert
function the ring $A$, the
monomials $e_{v,i,j}$, the multiindices $i$, $i'$, $i''$, triples $(v,i,j)$,
$(v,i',j')$, the
mo\-du\-le  $I$ and so on by the ring ${^h\!A}$, monomials $e_{v,i_0,i,j}$,
the pairs $(i_0,i)$, $(i'_0,i')$, $(i''_0,i'')$ (they are used without
parentheses), quadruples $(v,i_0,i,j)$,
$(v,i'_0,i',j')$, the
homogenization ${^h\!I}$ and so on respectively. Thus, we get the
definitions of $o(f)$, $\mathop{\rm Hdt}\nolimits(f)$ for
$f\in{^h\!A^l}$, new conditions (a) and (b) which define admissible
linear order on the monomials of ${^h\!A^l}$, new conditions
1)--5), the definitions of the Janet basis and reduced Janet basis
of ${^h\!I}$.  For example, the new conditions (a) and (b) are
\begin{enumerate} \renewcommand{\labelenumi}{(\alph{enumi})}
\item for all indices $i_0,i'_0$, all
multiindices $i,j,i',j'$ for all $1\le v\le l$
if $i_0\le i'_0$, $i_1\le i'_1,\ldots , i_n\le i'_n$ and $j_1\le j'_1,\ldots
, j_n\le j'_n$
then $(v,i_0,i,j)\le(v,i'_0,i',j')$.
\item for all indices $i_0,i'_0,i''_0$,
all multiindices $i,j,i',j',i'',j''$ for all $1\le v,v'\le l$
if  $(v,i_0,i,j)<(v',i'_0,i',j')$ then
$(v,i_0+i''_0,i+i'',j+j'')<(v',i'_0+i''_0,i'+i'',j'+j'')$.
\end{enumerate}
The Janet basis of ${^h\!I}$ is homogeneous if and
only if it consists of homogeneous elements from ${^h\!A^l}$.

Let $<$ be an admissible linear order on the monomials from $A^l$, or which
is the same, on the
triples $(v,i,j)$, see \srf{1}. So $<$ satisfies conditions (a) and (b).
Let us define the linear order on the monomials $e_{v,i_0,i,j}$ or, which
is the same, on the quadruples $(v,i_0,i,j)$. This linear order is induced
by  $<$ on the triples $(v,i,j)$ and will be denoted again by $<$.
Namely, for two quadruples $(v,i_0,i,j)$ and
$(v',i'_0,i',j')$ put $(v,i_0,i,j)<(v',i'_0,i',j')$
if and only if $(v,i,j)<(v',i',j')$, or $(v,i,j)=(v',i',j')$ but $i_0<i'_0$.
Notice that this induced linear order satisfies conditions (a) and (b)
(in the new sense).

\br{3} If $f_1,\ldots ,f_m$ is a Janet basis of $I$
(respectively homogeneous Janet basis of ${^h\!I}$)
satisfying
1)--4)
then there are the unique
$c_{\alpha,\beta}\in A$ (respectively $c_{\alpha,\beta}\in {^h\!A}$),
$1\le \alpha<\beta\le m$, such that
\[
f_\alpha+\sum_{\alpha<\beta\le m}c_{\alpha,\beta}f_\beta,\quad 1\le
\alpha\le m,
\]
is a reduced Janet basis of $I$ (respectively reduced homogeneous Janet basis
of
${^h\!I}$), cf. \cite{Cox}.
\end{rems}
%metka r3 \\

\bl{15} Let $f_1,\ldots , f_m$ be a
(reduced) Janet basis of $I$ with respect to the linear order $<$.
Then ${^h\!f}_1,\ldots ,{^h\!f}_m$ is a (reduced) homogeneous Janet basis of
the mo\-du\-le ${^h\!I}$ with respect to the induced linear order $<$.
Conversely, let $g_1,\ldots ,g_m$ be a (reduced) homogeneous Janet basis of
the mo\-du\-le ${^h\!I}$ with respect to the induced linear order $<$.
Then ${^a\!g}_1,\ldots ,{^a\!g}_m$ is a (reduced) Janet basis of $I$
with respect to the linear order $<$.
\end{lems}
%metka l15 \\
\noindent{\bf PROOF}\quad This follows immediately from the definitions.

\medskip\noindent Let $f\in{^h\!A^l}$ and the mo\-du\-le ${^h\!I}$ be as
above.
Then there is the unique element $g\in{^h\!A^l}$ such that
\[
g=\sum_{v,i_0,i,j}g_{v,i_0,i,j}e_{v,i_0,i,j},\quad g_{v,i_0,i,j}\in F,
\]
$f-g\in{^h\!I}$ and if $g_{v,i_0,i,j}\ne 0$ then
$e_{v,i_0,i,j}\not\in\mathop{\rm Hdt}\nolimits({^h\!I})$. The element
$g$ is called the normal form of $f$ with respect to the
mo\-du\-le ${^h\!I}$. We shall denote $g=\mathop{\rm
nf}\nolimits({^h\!I}, f)$. Obviously  $\mathop{\rm
nf}\nolimits({^h\!I},({^h\!A^l})_m)
\subset({^h\!A^l})_m$ is a linear subspace.

Let ${^c\!A}=F[X_0,\ldots , X_n,D_1,\ldots , D_n]$ be the polynomial ring in
the variables $X_0,\ldots , X_n,D_1,\ldots , D_n$.
Each monomial $e_{v,i_0,i,j}$ can be considered also as an element of
${^c\!A^l}$.
Denote by ${^c\!I}\subset{^c\!A^l}$ the graded submodule of ${^c\!A^l}$
generated by all the monomials $e_{v,i_0,i,j}$ such that there is
$0\ne f\in{^h\!I}$ with $o(f)=(v,i_0,i,j)$. The Hilbert functions
\begin{eqnarray*}
&&H({^c\!I},m)=\dim_F\{(z_1,\ldots , z_l)\in{^c\!I} \, :\,\forall\,
i\,(\,\deg
z_i=m\quad\mbox{or}\quad z_i=0\,)\},\\
&&H({^c\!A^l}/{^c\!I},m)={m+2n \choose 2n}-H({^c\!I},m).
\end{eqnarray*}
Let us replace in the definition of the normal form above ${^h\!A},{^h\!I}$
by ${^c\!A},{^c\!I}$
respectively. Thus, for $f\in{^c\!A^l}$
we get the definition of the normal form $\mathop{\rm nf}\nolimits({^c\!I},
f)\in
{^c\!A^l}$, cf. \cite{Dube}. Obviously, $\mathop{\rm
nf}\nolimits({^c\!I},({^c\!A^l})_m)
\subset({^c\!A^l})_m$ is a linear subspace. Since the ideals ${^c\!I}$ and
$\mathop{\rm Hdt}\nolimits({^h\!I})$ are generated by the same
monomials we have $\dim\mathop{\rm
nf}\nolimits({^c\!I},({^c\!A^l})_m)=\dim\mathop{\rm
nf}\nolimits({^h\!I},({^h\!A^l})_m)$.
Hence the Hilbert functions
\[
H({^h\!A^l}/{^h\!I},m)=H({^c\!A^l}/{^c\!I},m),
\quad H({^h\!I},m)=H({^c\!I},m),\quad m\ge 0,
\]
coincide.  Therefore, see \srf{3},
\bq{99}
H(I,m)=H({^c\!I},m),\quad m\ge 0
\end{equation}
%metka 99 \\

\section{Bound on the kernel of a matrix
over the homogenized Weyl algebra}\label{s4}
%metka s4 \\

\bl{4} Let $k=l-1$ and $l\ge 1$ be integers. Let $b=(b_{i,j})_{1\le i\le k,\,
1\le j\le l}$ be a ma\-t\-rix where $b_{i,j}\in {^h\!A}$ are homogeneous
elements
for all $i,j$. Let $\deg b_{i,j}<d$, $d\ge 1$, for all $i,j$.
Assume that there are integers $d_j\ge 0$, $1\le i\le k$, and
$d'_i\ge 0$, $1\le j\le l$, such that
\bq{66}
\deg b_{i,j}=d_i-d'_j
\end{equation}
%metka 66 \\
for all nonzero $b_{i,j}$, and
additionally $\min_{1\le j\le l}\{d'_j\}=0$ (hence $d_i<d$, $d'_j<d$ for all
$i,j$),
$d\ge 1$.
Then there are homogeneous elements
$z_1,\ldots , z_l\in{^h\!A}$ such that
$(z_1,\ldots, z_l)\ne(0,\ldots , 0)$,
\bq{32}
\sum_{1\le j\le l}b_{i,j}z_j=0,\quad 1\le i\le l-1,
\end{equation}
%metka 32 \\
all nonzero $b_{i,j}z_j$ have the same degree depending only on $i$ and
\bq{19}
\deg z_j\le (2n+3)ld,\quad 1\le j\le l.
\end{equation}
%metka 19 \\
Besides that, if all $b_{i,j}$ do not depend
on $X_n$
(i.e., they can be represented as sums of monomials which do not contain
$X_n$)
then one can choose also $z_1,\ldots , z_l$ satisfying additionally
the same property. Finally,
dividing by an appropriate power of $X_0$ one can assume without loss of
generality that
$\min\{\mathop{\rm ord}\nolimits z_i\, :\, 1\le i\le l\}=0$.
\end{lems}
%metka l4 \\
\noindent{\bf PROOF}\quad We shall assume without loss of generality that
$l\ge 2$. At first suppose that that \(\deg b_{i,j}=\deg b\) for all
nonzero $b_{i,j}$.
Consider the linear mapping
\bq{17}
\begin{array}{l}
({^h\!A})_{m-\deg b}^l
\longrightarrow({^h\!A})_m^{l-1}, \\
\left(\,z_1,\ldots , z_l\,\right)\mapsto \left(\,
\sum_{1\le j\le l}b_{i,j}z_j\,\right)_{1\le i\le l-1} .
\end{array}
\end{equation}
%metka 17 \\
If
\bq{13}
l{m-\deg b+2n \choose 2n}>(l-1){m+2n \choose 2n}
\end{equation}
%metka 13 \\
then the kernel of  \rf{17} is nonzero.
But \rf{13} holds if
\bq{14}
\Bigl(\,1+\frac{\deg b}{m+2n-\deg b}\,\Bigr)
\Bigl(\,1+\frac{\deg b}{m+2n-1-\deg b}\,\Bigr)\ldots
\Bigl(\,1+\frac{\deg b}{m-\deg b}\,\Bigr)<\frac{l}{l-1} .
\end{equation}
%metka 14 \\
Further,  \rf{14} is true if $(1+\deg b/(m-\deg b))^{2n}<l/(l-1)$.
The last inequality follows from $m\ge(2n+1)\deg b/\log(l/(l-1))$.
Hence also from $m\ge(2n+1)l\deg b$.
Notice that $(2n+2)ld\ge 1+(2n+1)l\deg b$.
Thus, the existence of $z_1,\ldots, z_l$ is proved, and even more
all nonzero $b_{i,j}z_j$ have the same degree which does not depend on $i,j$.
Notice that in the considered case we prove a more strong inequality $\deg
z_j\le(2n+2)ld$ for all $1\le j\le l$.

Suppose that $a_1,\ldots , a_l$ do not depend on $X_n$. We represent
$z_i=\sum_{j}z_{i,j}X_n^j$,
$1\le i\le l$, where all $z_{i,j}$ do not on $X_n$.
Let $\alpha=\max_i\{\deg_{X_n}z_i\}$.
Obviously in this case
one can replace $(z_1,\ldots , z_l)$ by $(z_{1,\alpha},\ldots ,
z_{l,\alpha})$.

Let us return to general case of arbitrary $\deg b_{i,j}$.
We shall reduce it to the considered one. Namely,
multiplying the  $i$-th
equation of system \rf{32} to
$X_0^{\max_i\{d_i\}-d_i}$
we shall suppose without loss of generality that all $d_i$ are
equal. Let us
substitute $z_jX_0^{d'_j}$ for $z_j$ in \rf{32}. Now the degrees
of all the nonzero coefficients of the obtained system coincide.
Thus, we get the required reduction and estimation  \rf{19}.
The lemma is proved.

\br{1} \lrf{4} remains
true if one replaces in its statement condition
\rf{32} by
\bq{33}
\sum_{1\le j\le l}z_jb_{i,j}=0,\quad 1\le i\le l-1,
\end{equation}
%metka 33 \\
The proof is similar.
\end{rems}
%metka r1 \\

\br{5} Let the elements $b_{i,j}$ be from \lrf{4}. Notice that there are
integers
$\delta'_i\ge 0$, $1\le i\le k$, and $\delta_j\ge 0$,
$1\le j\le l$, such that
\[
\deg b_{i,j}=\delta_j-\delta'_i
\]
for all nonzero $b_{i,j}$, and $\min_{1\le i\le k}\{\delta'_i\}=0$.
Namely, $\delta'_i=-d_i+\max_{1\le i\le k}\{d_i\}$,
$\delta_j=-d'_j+\max_{1\le i\le k}\{d_i\}$.
\end{rems}
%metka r5 \\

\section{Transforming a ma\-t\-rix with coefficients from ${^h\!A}$ to the
trapezoidal form}\label{s6}
%metka s6 \\

\noindent Let $b$ be the ma\-t\-rix from \lrf{4} but
now $k,l$ are arbitrary. Hence \rf{66} holds.
Let  $b=(b_1,\ldots , b_l)$ where $b_1,\ldots , b_l\in {^h\!A^k}$ be the
columns of the ma\-t\-rix $b$ (notice that in \lrf{2} and \lrf{3} $b_i$ are
rows of size $l$; so now we change the notation).
By definition $b_1,\ldots , b_l$ are linearly independent over ${^h\!A}$ from
the right
(or just linearly independent if it will not lead to an ambiguity)  if and
only if
for all $z_1,\ldots , z_l\in{^h\!A}$ the equality $b_1z_1+\ldots + b_lz_l=0$
implies $z_1=\ldots = z_l=0$. By \rf{66}
in this definition one can consider only homogeneous $z_1,\ldots , z_l$.
For an arbitrary family
$b_1,\ldots , b_l$ from \lrf{4} (with arbitrary $k,l$)
one can choose a maximal linearly independent from the right subfamily
$b_{i_1},\ldots , b_{i_r}$
of $b_1,\ldots , b_l$.
It turns out that $r$ does not depend on the choice of a subfamily.
More precisely, we have the following lemma.

\bl{11} Let $c_j=\sum_{1\le i\le l}b_iz_{i,j}$,  $1\le j\le r_1$, where
$z_{i,j}\in {^h\!A}$ are homogeneous elements. Suppose
that there are integers $d''_j$,  $1\le j\le r_1$, such that
for all $i,j$ the degree
$\deg z_{i,j}=d'_i-d''_j$. Assume that
$c_j$,
$1\le j\le r_1$, are linearly independent over ${^h\!A}$ from the right. Then
$r_1\le r$, and if $r_1<r$ there are $c_{r_1+1},\ldots ,
c_r\in\{b_{i_1},\ldots, b_{i_r}
\}$ such that $c_j$, $1\le j\le r$, are linearly independent over
${^h\!A}$ from the right.
\end{lems}
%metka l11 \\
\noindent{\bf PROOF}\quad The proof
is similar to the case of vector spaces over a field and we leave it to the
reader.

\medskip\noindent We denote $r=\mathop{\rm rankr}\nolimits\{b_1,\ldots ,
b_l\}$ and call it
the rank from the right of
$b_1,\ldots , b_l$. In the similar way one can define rank from the left of
$b_1,\ldots , b_l$. Denote it by $\mathop{\rm rankl}\nolimits\{b_1,\ldots ,
b_l\}$. It is not difficult to
construct examples when $\mathop{\rm rankr}\nolimits\{b_1,\ldots ,
b_l\}$ $\ne\mathop{\rm rankl}\nolimits\{b_1,\ldots , b_l\}$.
The aim of this section is to prove the following result.

\bl{9} Let $b$ be the ma\-t\-rix with homogeneous coefficient from ${^h\!A}$
satisfying \rf{66}, see above.
Suppose that $\deg b_{i,j}$ $<d$ for all $i,j$.
Assume that $k\ge l\ge 1$.
Let $l_1=\mathop{\rm rankr}\nolimits\{b_1,\ldots, b_l\}$ and
$b_1,\ldots , b_{l_1}$ be linearly independent. Hence $0\le l_1\le l$.
Then there is a ma\-t\-rix $(z_{j,r})_{1\le j,r\le l_1}$ with
homogeneous entries
$z_{j,r}\in{^h\!A}$ and a square permutation ma\-t\-rix $\sigma$ of
size $k$ satisfying the following properties.
\begin{enumerate} \renewcommand{\labelenumi}{(\roman{enumi})}
\item All the nonzero elements $b_{i,j}z_{j,r}$
for $1\leq j\leq l$ have the same degree
depending only on $i,r$ and
\bq{56}
\deg z_{j,r}\le (2n+3)ld.
\end{equation}
%metka 56 \\
\item Set the ma\-t\-rix $e=(e_{i,j})_{1\le i\le k,\,1\le
j\le l_1}=\sigma bz$. Then
the ma\-t\-rix
\[
e=\left(\begin{array}{ll}
e'\\
e''
\end{array}\right),
\]
where  $e'=\mathop{\rm diag}\nolimits(e'_{1,1},\ldots , e'_{l_1,l_1})$ is a
diagonal
ma\-t\-rix with $l_1$ columns and each $e'_{j,j}$,
$1\le j\le l_1$, is nonzero.
\item $\mathop{\rm ord}\nolimits e_{i,j}\ge
\mathop{\rm ord}\nolimits e'_{j,j}$ for all
$1\le i\le k$, $1\le j\le l_1$.
\end{enumerate}
Besides that, if all $a_{i,j}$ (and hence all $b_{i,j}$)
do not depend
on $X_n$
(i.e., they can be represented as sums of monomials which do not contain
$X_n$)
then one can choose also $z_{j,r}$ satisfying additionally
the same property. Finally,
dividing by an appropriate power of $X_0$ one can assume without loss of
generality that
$\min\{\mathop{\rm ord}\nolimits z_{j,r}\, :\,
1\le j\le l_1\}=0$ for every $1\le r\le l_1$.
\end{lems}
%metka l9 \\
\noindent{\bf PROOF}\quad At first we shall show how to construct $z$ and $e$
such that (ii) and (iii) hold.
We shall use a kind of Gauss elimination and \lrf{4}.
Namely, we transform the matrix $e$.
At the beginning we put
\[
e=(e_1,\ldots , e_{l_1})=(b_1,\ldots ,b_{l_1}).
\]
We shall perform some
${^h\!A}$-linear transformations of
columns and permutations of rows of the ma\-t\-rix $e$ and replace each time
$e$ by the obtained ma\-t\-rix. These transformation do not change
the rank from the right of the family of columns of $e$.
At the end we get a matrix $e$ satisfying the
required properties (ii), (iii).

We have $\mathop{\rm rankr}\nolimits(e)=l_1$.
If $l_1=0$, i.e, $e$ is an empty ma\-t\-rix, then this is the end
of the construction: $z'$ is an empty ma\-t\-rix.
Suppose that $l_1>0$.
Let us choose indices $1\le i_0\le k$, $1\le j_0\le l_1$
such that $\mathop{\rm ord}\nolimits e_{i_0,j_0}=
\min_{1\le j\le l_1}\{\mathop{\rm ord}\nolimits e_j\}$.
Permuting rows and columns of $e$ we shall assume without loss of generality
that $(i_0,j_0)=(1,1)$.

By \lrf{4} we get elements $w_{i,1},w_{i,i}\in{^h\!A}$ of
degrees at most $(2n+3)2d$ such that $e_{1,1}w_{1,i}=
e_{1,i}w_{i,i}$, $1\le i\le l_1$, and $\mathop{\rm
ord}\nolimits w_{i,i}=0$ for every $1\le i\le l_1$.
Set $w'=(-w_{1,2},\ldots ,-w_{1,l_1})$, and $w''=\mathop{\rm
diag}\nolimits(w_{2,2},\ldots , w_{l_1,l_1})$
to be the diagonal ma\-t\-rix. Put
\[
w=\left(\begin{array}{ll}
1, &w' \\
0, &w''
\end{array}\right)
\]
to be the square ma\-t\-rix with $l_1$ rows.
We replace $e$ by $ew$. Now
\[
e=\left(\begin{array}{ll}
e_{1,1}, &0\\
E_{2,1}, &E_{2,2}
\end{array}\right),
\]
where $E_{2,2}$ has $l_1-1$ columns and
\bq{16}
\min_{1\le j\le l_1}\{\mathop{\rm ord}\nolimits b_j\}=
\mathop{\rm ord}\nolimits e_{1,1}=
\min_{1\le j\le l_1}\{\mathop{\rm ord}\nolimits e_j\}
\end{equation}
%metka 16 \\
(for the new ma\-t\-rix $e$).

Let us apply recursively the described construction to the ma\-t\-rix
$E_{2,2}$ in place of $e$.
So using only linear transformations of columns with indices $2,\ldots , l_1$
and permutation of rows with indices $2,\ldots , k$ we transform $e$ to the
form
\[
\sigma e\tau=\left(\begin{array}{ll}
e_{1,1}, &0 \\
E'_{2,1}, &E'_{2,2} \\
E''_{2,1} &E''_{2,2}
\end{array}\right), \quad
\tau=\left(\begin{array}{ll}
1, &0 \\
0, &\tau'
\end{array}\right)
\]
where $\sigma$ is a permutation ma\-t\-rix and $\tau'$ is a square ma\-t\-rix
with $l_1-1$ rows (it transforms $E_{2,2}$), the ma\-t\-rix
$E'_{2,2}=\mathop{\rm diag}\nolimits(e_{2,2},\ldots , e_{l_1,l_1})$ is
a diagonal ma\-t\-rix with
$l_1-1\ge 0$ columns, and all the elements $e_{2,2},\ldots , e_{l_1,l_1}\in
{^h\!A}$
are nonzero. We shall assume without loss of generality that $\sigma=1$ is
the
identity ma\-t\-rix. We replace $e$ by $e\tau$.
Conditions (ii) and (iii) hold for the obtained $e$ and, more
than that,
by (iii) applied recursively for $(E_{2,2},E'_{2,2},E''_{2,2})$
(in place of $(e,e',e'')$), and
\rf{16} the same equalities are satisfied
for the new obtained
ma\-t\-rix $e$.

Let $E'_{2,1}=(e_{2,1},\ldots , e_{l_1,1})^t$ where $t$ denotes
transposition.
By \lrf{4} there are nonzero
elements $v_{1,1},\ldots , v_{l_1,1}\in{^h\!A}$
of degrees at most
\bq{21}
(2n+3)(\max\{\deg e_{i,i}\, :\, 1\le i \le
l_1\}+1)l_1
\end{equation}
%metka 21 \\
such that $e_{i,1}v_{1,1}=e_{i,i}v_{i,1}$ and $\min\{\mathop{\rm
ord}\nolimits v_{1,1},\mathop{\rm ord}\nolimits v_{1,i}\}=0$ for
all $1\le i\le l_1-1$.
Let $v'=(-v_{2,1},\ldots ,-v_{l_1,1})^t$ and $v''$
be the identity
ma\-t\-rix of
size $l_1-1$.
Put
\[
v=\left(\begin{array}{ll}
v_{1,1}, &0 \\
v', &v''
\end{array}\right).
\]
Let us replace $e$ by $ev$. Put $z=w\tau v$, where
the ma\-t\-rix $z$ has $l_1$ columns. Recall that without loss of
generality  $\sigma=1$
is the identity permutation. We have $e=(b_1,\ldots , b_{l_1})z$.
These Gauss elimination transformations of $e$ do not change
the rank from the right of the family of columns of $e$.
It can be easily proved using the recursion on $l$, cf. \lrf{7} below.
Now the ma\-t\-rix $e$ satisfies required
conditions (ii), (iii) and $\sigma=1$.

Let us change the notation. Denote the obtained ma\-t\-rix $z$ by $z'$.
Let $z'=(z'_1,\ldots , z'_{l_1})$ where $z'_j$ is the $j$-th column of $z'$.
Our aim now is to prove the existence of the ma\-t\-rix $z$ satisfying
(i)--(iii).
By \lrf{4} for every $1\le r\le l_1$ there are homogeneous elements
$z_{j,r}\in {^h\!A}$, $1\le
j\le l$,
such that  $(z_{1,r},\ldots , z_{l,r})\ne(0,\ldots , 0)$,
\bq{37}
\sum_{1\le j\le l_1}b_{i,j}z_{j,r}=0\quad\mbox{for every}\quad 1\le i\le
l_1,\;i\ne r,
\end{equation}
%metka 37 \\
and estimations for degrees \rf{56} hold. Put the ma\-t\-rix
$z=(z_{j,r})_{1\le j,r\le l_1}$.
Let $z=(z_1,\ldots , z_{l_1})$ where $z_j$ is the $j$-th column of $z$.
Hence $z_j=(z_{1,r},\ldots , z_{l,r})^t$.

\bl{10} For every $1\le r\le l_1$ we have
\bq{35}
\sum_{1\le j\le l_1}b_{r,j}z_{j,r}\ne 0.
\end{equation}
%metka 35 \\
Further, for every $1\le r\le l_1$
there are nonzero homogeneous elements $g'_r,g_r\in{^h\!A}$ such that
$z'_rg'_r=z_rg_r$.
\end{lems}
%metka l10 \\
\noindent{\bf PROOF}\quad
Consider the matrix $(z',z_r)$ with  $l_1$ rows and
$l_1+1$ columns. By \lrf{4} there are homogeneous elements
$h_1,\ldots , h_{l_1+1}\in{^h\!A}$ (they depend on $r$)
such that $(h_1,\ldots , h_{l_1+1})\ne(0,\ldots , 0)$ and the following
property holds. Denote $h=(h_1,\ldots, h_{l_1+1})^t$, $h'=(h_1,\ldots,
h_{l_1})^t$.
Then
\bq{36}
z'h'+z_rh_{l_1+1}=0
\end{equation}
%metka 36 \\
(we don't need at present any estimation on
degrees
from \lrf{4}; only the existence of $h$). Denote by $b''$ the submatrix
consisting of
the first $l_1$ rows of the ma\-t\-rix $(b_1,\ldots , b_{l_1})$.
Multiplying  \rf{36} to $b''$
from the left we get
\bq{38}
b''z'h'+b''z_rh_{l_1+1}=0.
\end{equation}
%metka 38 \\
But $b''z'$ is a diagonal ma\-t\-rix with
nonzero elements on the
diagonal, see (ii) (for $z'$ in place of $z$). Hence by  \rf{37} and \rf{38}
$h_j=0$ for every $j\ne r$. Now  $h\ne(0,\ldots , 0)^t$
implies $h_r\ne 0$ and $h_{l_1+1}\ne 0$. Therefore, \rf{35} holds.
Put $g'_r=h_r$ and $g_r=h_{l_1+1}$. We have $z'_rg'_r=z_rg_r$ by  \rf{38}.
The lemma is proved.

\medskip\noindent  Let us return to the proof of \lrf{9}. Now (i)--(iii)
are satisfied by \lrf{10}.
The last assertions of \lrf{9} are proved similarly to the ones of  \lrf{4}.
\lrf{9} is proved.

\section{An algorithm for solving linear systems with
coefficients from ${^h\!A}$.}\label{s5}
%metka s5 \\

\noindent Let $u=(u_1,\ldots, u_l)^t\in{^h\!A^l}$. Let all nonzero $u_j$ be
homogeneous elements of the degree $-d'_j+\rho$ for an integer $\rho$.
Suppose that $-d'_j+\rho<d'$ for an integer $d'>1$.
Let $b=(b_{i,j})_{1\le i\le k,\,1\le j\le l}$ be the ma\-t\-rix
with  $k$ rows and $l$ columns from the
statement of \lrf{9} (but now $k$ and $l$ are arbitrary). So $\deg
b_{i,j}=d_i-d'_j<d$ for all $i,j$.
Let $Z=(Z_1,\ldots , Z_k)$ be unknowns.
Consider the linear system
\bq{22}
\sum_{1\le i\le k}Z_i b_{i,j}=u_j,\quad 1\le j\le l,
\end{equation}
%metka 22 \\
or, which is the same,
\[
Zb=u.
\]
Denote
\bq{57}
\mathop{\rm ord}\nolimits u=\min_{1\le i\le k}
\{\mathop{\rm ord}\nolimits u_i\}.
\end{equation}
%metka 57 \\
The similar notations will be
used for other
vectors and ma\-t\-ri\-ces.
In this section we shall describe an algorithm for solving linear systems
over ${^h\!A}$ and prove the following theorem.

\bthm{2} Suppose that system \rf{22} has a solution over ${^h\!A}$.
One can represent the set
of all solutions of \rf{22} over ${^h\!A}$ in the form
\[
J+z^*,
\]
where $J\subset {^h\!A^l}$ is a ${^h\!A}$-submodule of all the solutions of
the homogeneous system corresponding to \rf{22}
(i.e., system  \rf{22} with all $u_j=0$)
and $z^*$ is a particular solution of \rf{22}.
Moreover, the following assertions hold.
\begin{enumerate} \renewcommand{\labelenumi}{(\Alph{enumi})}
\item One can choose $z^*$ such that
$\mathop{\rm ord}\nolimits z^*\ge
\mathop{\rm ord}\nolimits u-\nu$,
where $\nu\ge 0$ is an integer
bounded from above by $(dl)^{2^{O(n)}}$ (and depends
only on $d$ and $l$).
The degree $\deg z^*$ is bounded from above by $d'(dl)^{2^{O(n)}}$.
\item There exists a system of generators of $J$ of degrees
bounded
from above by $(dl)^{2^{O(n)}}$. The number of elements of this system of
generators is bounded from above by $k(dl)^{2^{O(n)}}$.
\end{enumerate}
Besides that, if all $b_{i,j}$ and $u_j$ do not depend
on $X_n$
(i.e., they can be represented as sums of monomials which do not contain
$X_n$)
then $z^*$ and all the generators of the mo\-du\-le $J$ also satisfy this
property.
\end{thms}
%metka t2 \\

\medskip\noindent{\bf PROOF}\quad
Let $l_1=\mathop{\rm rankr}\nolimits(b_1,\ldots , b_l)$. Permuting equations
of  \rf{22} we shall
assume without loss of generality
that $(b_1,\ldots , b_{l_1})$ are linearly independent from the right over
${^h\!A}$.
Let $\sigma,z,e,e',e''$ be the ma\-t\-ri\-ces from  \lrf{9}.
Similarly to the proof of \lrf{9} we shall assume without loss of generality
that $\sigma=1$.
Denote by $b'$ the submatrix of $b$ consisting of the first $l_1$ columns of
$b$,  i.e.,
$b'=(b_1,\ldots , b_{l_1})$.
By \lrf{4} there are nonzero elements $q_{1,1},\ldots ,$ $q_{l_1,l_1}$ of
degrees at most
\rf{21} such that $e_{1,1}q_{1,1}=e_{i,i}q_{i,i}$ and $\min\{\mathop{\rm
ord}\nolimits q_{1,1},\mathop{\rm ord}\nolimits q_{i,i}\}=0$
for all $2\le i\le l_1$. Set $q=\mathop{\rm diag}\nolimits(q_{1,1},\ldots ,
q_{l_1,l_1})$ to be the diagonal ma\-t\-rix.
Let $\nu_0=\mathop{\rm ord}\nolimits e_{1,1}q_{1,1}$. Then by \lrf{9}
(iii)
$\mathop{\rm ord}\nolimits (b'zq)\ge\nu_0$.
Let $X_0^{\nu_0}\delta=b'zq$. Then $\delta$ is a ma\-t\-rix with coefficients
from ${^h\!A}$ and
\[
\delta=
\left(
\begin{array}{c}
\delta' \\
\delta''
\end{array}
\right),
\]
where $\delta'=\mathop{\rm diag}\nolimits(\delta_{1,1},\ldots
,\delta_{l_1,l_1})$ is a diagonal ma\-t\-rix with homogeneous coefficients
from ${^h\!A}$
and all the elements on the diagonal are nonzero and equal,
i.e., $\delta_{j,j}=\delta_{1,1}$ for every
$1\le j\le l_1$. Besides that, $\mathop{\rm
ord}\nolimits \delta_{1,1}=0$.  Further,
$\delta''=(\delta_{i,j})_{l_1+1\le i\le k,\,
1\le j\le l_1}$.
We have $\mathop{\rm ord}\nolimits (uzq)\ge\nu_0$,
since, otherwise, system \rf{22} does not
have a solution. Obviously $\mathop{\rm ord}\nolimits u\le\mathop{\rm
ord}\nolimits(uzq)$.
Denote $u'=(u'_0,\ldots , u'_l)^t=X_0^{-\nu_0}uzq\in {^h\!A^l}$.
Hence $\mathop{\rm ord}\nolimits u'\ge\mathop{\rm ord}\nolimits (u)-\nu_0$.
Consider the linear system
\bq{23}
Z\delta=u'.
\end{equation}
%metka 23 \\

\bl{7} Suppose that system \rf{22} has a solution over ${^h\!A}$.
Then linear system  \rf{23} is equivalent to  \rf{22}, i.e., the sets of
solutions of systems \rf{23} and  \rf{22} over ${^h\!A}$ coincide.
\end{lems}
%metka l7 \\
\noindent{\bf PROOF}\quad The system $Zb'z=uz$ is equivalent to \rf{22} by
\lrf{11}.
System  \rf{23} is equivalent to $Zb'z=uz$ since the ring ${^h\!A}$
does not have zero--divisors.
The lemma is proved.

\br{4} Since $\mathop{\rm rankr}\nolimits(b_1,\ldots , b_l)=l_1$ and
by \lrf{9} for every $l_1+1\le j\le l$ there are homogeneous
$z_{j,j},z_{1,j},\ldots , z_{l_1,j}\in{^h\!A}$ such that  $z_{j,j}\ne 0$ and
$b_jz_{j,j}+\sum_{1\le r\le l_1}b_rz_{r,j}=0$
and all $\deg z_{j,j}$, $\deg z_{r,j}$
are bounded from above by $(2n+3)(l_1+1)d$. Put
$u'_j=u_jz_{j,j}+\sum_{1\le r\le l_1}u_rz_{r,j}$, $l_1+1\le j\le l$.
Then system  \rf{22} has a solution if and only if system \rf{23} has a
solution and $u'_j=0$
for all $l_1+1\le j\le l$. This follows from \lrf{7} and \lrf{11}.
But in what follows for our aims it is sufficient to use only \lrf{7}.
\end{rems}
%metka r4 \\

\br{2} Assume that  $\deg_{X_n}b_{i,j}\le 0$ for all $i,j$, i.e., the
elements of the ma\-t\-rix $b$ do not depend on $X_n$. Then by \lrf{4}
and the described construction all the elements of the ma\-t\-ri\-ces
$b,z,q, \delta, \delta', \delta''$ also
do not depend on $X_n$.
\end{rems}
%metka r2 \\

By \lrf{4} and \rrf{1} for every $l_1+1\le j\le k$ there are  homogeneous
elements $g_{j,j},g_{j,i}\in {^h\!A} $, $1\le i\le l_1$, such that
\[
g_{j,j}\delta_{j,i}=g_{j,i}\delta_{1,1},\quad 1\le i\le l_1,
\]
all the degrees $\deg g_{j,j}$, $\deg g_{j,i}$, $1\le i\le l_1$, are bounded
from above by
\[
(2n+3)(l_1+1)(\max\{\deg \delta_{j,i}\, :\, 1\le i \le
k\}+1)
\]
and  $\min_{1\le i\le l_1}\{\mathop{\rm ord}\nolimits g_{j,j},
\mathop{\rm ord}\nolimits g_{j,i}\}=0$.
Hence $\mathop{\rm ord}\nolimits g_{j,j}=0$ for every $l_1+1\le j\le k$
since $\mathop{\rm ord}\nolimits \delta_{1,1}=0$.

Denote $h=\delta_{1,1}g_{l_1+1,l_1+1}g_{l_1+2,l_1+2}\ldots g_{k,k}$.
So $h\in {^h\!A}$ is a nonzero homogeneous element and $\mathop{\rm
ord}\nolimits h=0$.
Set $\varepsilon=\deg h$. We need an analog of the Noether
normalization theorem
from commutative al\-ge\-b\-ra, cf. also Lemma~3.1 \cite{G05}.

\bl{6} There is a linear automorphism of the
al\-ge\-b\-ra ${^h\!A}$
\begin{eqnarray*}
&&\alpha\, :\,{^h\!A}\rightarrow{^h\!A},\quad
\alpha(X_i)=\sum_{1\le j\le n}(\alpha_{1,i,j}X_j+\alpha_{2,i,j}D_j), \\
&&\alpha(D_i)=\sum_{1\le j\le
n}(\alpha_{3,i,j}X_j+\alpha_{4,i,j}D_j),\quad\alpha(X_0)=X_0,\quad 1\le
i\le n,
\end{eqnarray*}
such that all $\alpha_{s,i,j}\in F$, $\deg_{D_n}\alpha(h)=\varepsilon$.
If $\deg_{X_n}h=0$ then one can choose additionally $\alpha(X_n)=X_n$,
all
$\alpha_{1,n,j}=0$ for $1\le j\le n-1$ and $\alpha_{3,n,j}=0$ for $1\le j\le
n$.
\end{lems}
%metka l6 \\
\noindent{\bf PROOF}\quad Recall that $\mathop{\rm ord}\nolimits h=0$.
Hence at first it is not difficult to construct a linear
automorphism $\beta$ such that  $\beta(X_0)=X_0$,
\bq{61}
\beta(X_i)=\beta_{1,i}X_i+\beta_{2,i}D_i,\quad
\beta(D_i)=\beta_{3,i}X_i+\beta_{4,i}D_i,\quad 1\le i\le n,
\end{equation}
%metka 61 \\
and $\beta(h)$
contains a monomial $a_{i_1,\ldots , i_n}D_1^{i_1},\ldots , D_n^{i_n}$ with
$a_{i_1,\ldots , i_n}\ne 0$ and
$i_1+\ldots + i_n=\varepsilon$, i.e., $\varepsilon=\deg_{D_1,\ldots ,
D_n}\beta(h)$. After that one can find an automorphism
$\gamma$ such that $\gamma(X_0)=X_0$,
\bq{62}
\gamma(X_i)=\sum_{1\le j\le n}\gamma_{1,i,j}X_j,\quad
\gamma(D_i)=\sum_{1\le j\le n}\gamma_{4,i,j}D_j,\quad 1\le i\le n,
\end{equation}
%metka 62 \\
and $(\gamma\circ\beta)(h)$
contains a monomial  $aD_n^\varepsilon$ with a coefficient $0\ne a\in F$. Put
$\alpha=\gamma\circ\beta$. We leave to prove the last assertion to the
reader.
The lemma is proved.

\medskip\noindent We apply the automorphism $\alpha$.
In what follows to simplify the notation we shall suppose without loss of
generality that $\alpha=1$. So  $h$ contains a monomial
$aD_n^\varepsilon$ with a coefficient $0\ne a\in F$, where $\varepsilon=\deg
h$.
It follows from here that
\bq{24}
\deg_{D_n}\delta_{1,1}=\deg \delta_{1,1},\quad \deg_{D_n}g_{j,j}=
\deg g_{j,j},\,l_1+1\le j\le k.
\end{equation}
%metka 24 \\
Let $z=(z_1,\ldots , z_k)\in{^h\!A^k}$ be a solution of \rf{23}.
Then  \rf{24} implies that one can uniquely represent
\bq{26}
z_j=z'_jg_{j,j}+\sum_{0\le s<\deg g_{j,j}}z_{j,s}D_n^s,\quad l_1+1\le j\le k,
\end{equation}
%metka 26 \\
where $z'_j,z_{j,s}\in{^h\!A}$, the degrees $\deg_{D_n}z_{j,s}\le 0$ for all
$l_1+1\le j\le k$, $0\le
s<\deg_{D_1}g_{j,j}$.
Again by  \rf{24} one can uniquely represent
\[
u'_i=u''_i\delta_{1,1}+\sum_{0\le s<\deg \delta_{1,1}}u'_{i,s}D_n^s,\quad
1\le i\le l,
\]
where $u''_i,u'_{i,s}\in{^h\!A}$,
the degrees $\deg_{D_n}u'_{i,s}\le 0$ for all  $1\le i\le l$, $0\le
s<\deg_{D_1}g_{j,j}$.
Finally, by  \rf{24} for all $l_1+1\le j\le k$, $1\le i\le l_1$, $0\le
r<\deg_{D_1}g_{j,j}$, one can uniquely represent
\[
D_n^r\delta_{j,i}=\delta_{j,r,i}\delta_{1,1}+\sum_{0\le
r<\deg \delta_{1,1}}\delta_{j,r,i,s}D_n^s,
\]
where $\delta_{j,r,i},\delta_{j,r,i,s}\in{^h\!A}$,
the degrees $\deg_{D_n}\delta_{j,r,i,s}\le 0$ for all considered $j,r,i,s$.
Put
\begin{eqnarray*}
&&{\mathcal I}=\left\{\,(j,r)\, :\, l_1+1\le j\le k\,\&\,0\le r< \deg
g_{j,j}\,\right\}, \\
&&{\mathcal J}=\left\{\,(i,s)\,:\,1\le i\le l_1\,\&\,1\le s< \deg
\delta_{1,1}\,\right\}.
\end{eqnarray*}

Therefore,
\begin{eqnarray}
&&z_i=-\sum_{l_1+1\le j\le k}z'_jg_{j,i}-\sum_{(j,r)\in{\mathcal
I}}z_{j,r}\delta_{j,r,i}+u''_i,\quad 1\le i\le l_1, \label{025} \\
&&\sum_{(j,r)\in{\mathcal I}}z_{j,r}\delta_{j,r,i,s}=u'_{i,s},\quad
(i,s)\in{\mathcal J}. \label{026}
\end{eqnarray}
%metka 025,026 \\
Let us introduce new unknowns $Z_{j,r}$, $(j,r)\in{\mathcal I}$. By
 \rf{26}--\rf{026} system \rf{22} is reduced to the linear system
\bq{27}
\sum_{(j,r)\in{\mathcal I}}Z_{j,r}\delta_{j,r,i,s}=u'_{i,s},\quad
(i,s)\in{\mathcal J}.
\end{equation}
%metka 27 \\
More precisely, any solution of system \rf{22} is given by  \rf{26}, \rf{025}
where $z'_j\in{^h\!A}$
are arbitrary and $z_{j,r}$ is a solution of system  \rf{026} over ${^h\!A}$
(we
underline that here this solution $z_{j,r}$ may depend on $D_n$ although one
can restrict oneself by solutions $z_{j,r}$ which do not depend on $D_n$).
Note that all $\delta_{j,r,i,s}$ and $u'_{i,s}$
are homogeneous elements of ${^h\!A}$ and
there are integers $d_{j,r}$, $(j,r)\in{\mathcal I}$, $d'_{i,s}$,
$(i,s)\in{\mathcal J}$,  $\widetilde{\rho}$
such that $\deg\delta_{j,r,i,s}=d_{j,r}-d'_{i,s}$
and $\deg u'_{i,s}=-d'_{i,s}+\widetilde{\rho}$ for all $(j,r)\in{\mathcal
I}$, $(i,s)\in{\mathcal J}$.
This follows immediately from the described construction.

Now all the coefficients of system  \rf{27} do not depend on $D_n$.
As we have proved if the coefficients of  \rf{22} do not depend on $X_n$ then
the coefficients of  \rf{27} also do not depend on $X_n$, and hence in the
last case
they do not depend on $X_n,D_n$.

If the coefficients of \rf{27} depend on $X_n$ we perform an automorphism
$X_n\mapsto D_n$ $D_n\mapsto-X_n$, $X_i\mapsto X_i$, $D_i\mapsto D_i$, $1\le
i\le n-1$.
Now the coefficients of system  \rf{27} do not depend on $X_n$ (but depend
on $D_n$). After that we apply our construction recursively to system
\rf{27}.

The final step of the recursion is $n=0$ (although in the statement of
theorem $n\ge 1$, see \srf{1}; we are interested only in Weyl
al\-ge\-b\-ras). In this case ${\mathcal I}={\mathcal J}=\emptyset$.
Hence using \rf{025} for $n=0$ we get the required $z^*$ and $J$ for $n=0$.

Thus, by the recursive assumption we get a particular solution
$Z_{j,r}=z_{j,r}^*$, $
(j,r)\in{\mathcal I}$,
of system  \rf{27}, an integer $\nu_1$ (in place of $\nu$ from assertion
(A)) such that
\bq{105}
\min_{(j,r)\in{\mathcal I}}\{\mathop{\rm ord}\nolimits z_{j,r}^*\}\ge
\min_{(i,s)\in{\mathcal J}}\{\mathop{\rm ord}\nolimits u'_{i,s}\}-\nu_1,
\end{equation}
%metka 105 \\
and a system of generators
\bq{39}
\left(\,z_{\alpha,j,r}\,\right)_{(j,r)\in{\mathcal I}},\quad 1\le\alpha \le
\beta,
\end{equation}
%metka 39 \\
of the mo\-du\-le $J'$ of solutions of the homogeneous system
corresponding to \rf{27}. Notice that
if the coefficients of  \rf{22} do not depend on $X_n$ then
$J'$ is a mo\-du\-le over the homogenization $F[X_0,
X_1,\ldots,X_{n-1},
D_1,\ldots,D_{n-1}]$ of the Weyl al\-ge\-b\-ra of $X_1,\ldots,X_{n-1},
D_1,\ldots,D_{n-1}$.
But obviously in the last case  \rf{39} gives also a system of generators of
the ${^h\!A}$-mo\-du\-le $J''={^h\!A}J'$ of solutions of the homogeneous
system corresponding to \rf{27}.
Put
\begin{eqnarray*}
&&z^*_i=-\sum_{(j,r)\in{\mathcal I}}z^*_{j,r}\delta_{j,r,i}+u''_i,\quad 1\le
i\le l_1,\\
&&z^*_j=\sum_{0\le s<\deg g_{j,j}}z^*_{j,s}D_n^s,\quad l_1+1\le j\le k,\\
&&z^*=(z^*_1,\ldots , z^*_k).
\end{eqnarray*}
Then $z^*$ is a particular solution of  \rf{22}.
Put
\begin{eqnarray*}
&&z_{\alpha,i}=-\sum_{(j,r)\in{\mathcal
I}}z_{\alpha,j,r}\delta_{j,r,i},\quad 1\le i\le l_1,\;1\le\alpha\le
\beta,\\
&&z_{\alpha,j}=\sum_{0\le s<\deg g_{j,j}}z_{\alpha,j,s}D_n^s,\quad
l_1+1\le j\le k,\;1\le\alpha\le \beta, \\
&&z_{\beta-l_1+j,i}=0,\quad l_1+1\le i,j\le k,\; i\ne j,\\
&&z_{\beta-l_1+j,j}=g_{j,j},\quad l_1+1\le j\le k, \\
&&z_{\beta-l_1+j,i}=-g_{j,i},\quad 1\le i\le l_1,\; l_1+1\le j\le k.
\end{eqnarray*}
Then  $J=\sum_{1\le\alpha\le \beta+k-l_1}{^h\!A}(z_{\alpha,1},\ldots ,
z_{\alpha,k})$. Hence $(z_{\alpha,1},\ldots,
z_{\alpha,k})$,  $1\le\alpha\le\beta+k-l_1$, is a system of generators of the
mo\-du\-le $J$.
By  \rf{105} and the definitions of $u'$, $u''_i$ and $u'_{i,s}$ we have
$\mathop{\rm ord}\nolimits z^*\ge\mathop{\rm ord}\nolimits u-\nu_0-\nu_1$.
Put $\nu=\nu_0+\nu_1$.

\bl{8} All the degrees  $\deg\delta_{j,i}$, $\deg g_{j,i}$,
$\deg\delta_{j,r,i}$,
$\deg \delta_{j,r,i,s}$ and $\nu$, see above, are bounded from above by
$(nld)^{O(1)}$, the degrees $\deg u'_i$ are bounded from above
$d'+(nld)^{O(1)}$,
the degrees $\deg u''_i$, $\deg u'_{i,s}$ are bounded from above by
$d'(nld)^{O(1)}$.
Further, all $\mathop{\rm ord}\nolimits u''_i$, $\mathop{\rm
ord}\nolimits u'_{i,s}$ are bounded from below by
$\mathop{\rm ord}\nolimits u-\nu$.
Finally, in system  \rf{27} the number of equations  $\#{\mathcal J}$
is bounded from
above by $(nld)^{O(1)}$ and the number
of unknowns $\#{\mathcal I}$ is bounded from
above by $k(nld)^{O(1)}$.
\end{lems}
%metka l8 \\
\noindent{\bf PROOF}\quad This follows immediately from the described
construction.

\medskip\noindent Let us return to the proof of \trf{2}. Applying \lrf{8}
and recursively assertions (A) and (B) for
the formulas giving  $z^*$ and $J$
we get (A) and (B) from the theorem. The last
assertion (related  to the case when all $b_{i,j}$ and $u_j$ do not depend
on $D_n$)
has been already proved.
The theorem is proved.

\section{Proof of \trf{1} for Weyl algebra}\label{s9}
%metka s9 \\

\noindent Let $a$ be the ma\-t\-rix from \srf{1}.
We shall suppose without loss of generality that the vectors
$(a_{i,1},\ldots , a_{i,l})$, $1\le i\le k$,
are linearly independent over the field $F$. We have $\deg a_{i,j}<d$. This
implies
$k\le l{d+2n\choose 2n}$.

Put the ma\-t\-rix $b={^h\!a}$.
Let us define the graded submodules of ${^h\!I}$
\begin{eqnarray*}
&&J_0={^h\!A}(b_{1,1},\ldots , b_{1,l})+\ldots +{^h\!A}(b_{k,1},\ldots ,
b_{k,l}),  \\
&&J_\nu=J_0:(X_0^\nu)=\{z\in {^h\!A^l}\, :\, zX_0^\nu\in J_0
\},\quad
\nu\ge 1.
\end{eqnarray*}
We have the exact sequence of graded ${^h\!A}$-mo\-du\-les
\[
{^h\!A^k}\rightarrow J_0\rightarrow 0.
\]
Further, $J_\nu\subset J_{\nu+1}\subset{^h\!I}$ for every $\nu\ge
0$ and
${^h\!I}=\bigcup_{\nu\ge 0}J_\nu$.
Since ${^h\!A}$ is Noetherian there is $N\ge 0$ such that
${^h\!I}=J_N$. So to construct a system of generators of ${^h\!I}$ it is
sufficient to compute
the least $N$ such that ${^h\!I}=J_N$ and to find a system of generators of
$J_N$.

\bl{5} ${^h\!I}=J_N$ for some $N$ bounded from above by $(dl)^{2^{O(n)}}$.
There is a system of generators $b_1,\ldots , b_s$ of the mo\-du\-le $J_N$
such that $s$ and
all the degrees $\deg b_v$,
$1\le v\le s$, are bounded from above by $(dl)^{2^{O(n)}}$.
\end{lems}
%metka l5 \\
\medskip\noindent{\bf PROOF}\quad
Let us show that the mo\-du\-le $J_{N+1}\subset J_N$ for
$N\ge\nu$. Let $u\in J_{N+1}$. Consider system \rf{22}.
By assertion (A) of \trf{2} there is a particular solution $z^*$ of  \rf{22}
such that $\mathop{\rm ord}\nolimits z^*\ge 1$. Hence $u\in X_0J_N\subset
J_N$. The required assertion is proved. Hence ${^h\!I}=J_\nu$.

Let us replace in \rf{22} $(u_1,\ldots , u_l)$ by $(U_1X_0^\nu,\ldots ,
U_lX_0^\nu)$, where
$U_1,\ldots, U_l$ are new unknowns. Then applying (B) from \trf{2} to this
new
homogeneous linear system with respect to all unknowns $U_1,\ldots , U_l$,
$Z_1,\ldots ,Z_k$
we get the required estimations for the number of generators of $J_\nu$ and
the degrees of these generators. The lemma is proved.

\bco{1} Let $(a_{i,1},\ldots , a_{i,l})$, $1\le i\le l$, be from the
beginning of the section and the integer $N$ be from \lrf{15}. Then for
every integer $m\ge 0$ the $F$--linear
space
\bq{106}
A_{m+N}(a_{1,1},\ldots , a_{1,l})+\ldots + A_{m+N}(a_{k,1},\ldots ,
a_{k,l})\supset I_m.
\end{equation}
%metka 106 \\
\end{coros}
%metka c1 \\
\noindent{\bf PROOF}\quad By \lrf{15} we have $(J_0)_{m+N}\supset
X_0^N(J_N)_m=X_0^N({^h\!I})_m$. Taking the
affine parts we get  \rf{106}. The corollary is proved.

\medskip\noindent Now everything is ready for the proof of \trf{1}.
By \lrf{5} and \lrf{2} there is a system of
generators of the mo\-du\-le $\mathop{\rm gr}\nolimits(I)$ with
degrees bounded from above by $(dl)^{2^{O(n)}}$.
By \lrf{14} from Appendix~1 the Hilbert function $H(\mathop{\rm
gr}\nolimits(I),m)$
is stable for $m\ge(dl)^{2^{O(n)}}$. By \rf{65} \srf{2} the Hilbert function
$H(I,m)$ is stable for all $m\ge(dl)^{2^{O(n)}}$.

Consider the linear order $<$ on the monomials from ${^h\!A^l}$
which is induced by the linear order $<$ on the monomials from $A^l$, see
\srf{10}.
Then the monomial submodule ${^c\!I}\subset{^c\!A^l}$ is defined, see
\srf{10}, where
${^c\!A}=F[X_0,\ldots , X_n,D_1,\ldots ,$ $D_n]$ is the polynomial ring.
By \rf{99} \srf{10} the Hilbert function $H({^c\!I},m)$
is stable for all $m\ge(dl)^{2^{O(n)}}$. Hence all the coefficients of the
Hilbert polynomial
of ${^c\!I}$ are bounded from above $(dl)^{2^{O(n)}}$.
Therefore, according to \rf{16} the mo\-du\-le ${^c\!I}$ has a system of
generators with
degrees $(dl)^{2^{O(n)}}$. This means, see \srf{10}, that the mo\-du\-le
$\mathop{\rm Hdt}\nolimits({^h\!I})$ has a system of generators with
degrees $(dl)^{2^{O(n)}}$. Therefore, the degrees of all the elements of the
Janet basis of ${^h\!I}$ with respect to the induced linear order $<$ are
bounded from above by $(dl)^{2^{O(n)}}$.
Hence by \lrf{15} \srf{10} the same is true for the Janet basis of the
mo\-du\-le $I$ with respect to the linear order $<$ on the monomials from
$A^l$. \trf{1} is proved for Weyl algebra.

\section{The case of al\-ge\-b\-ra of differential operators}\label{s8}
%metka s8 \\

Denote by $B=F(X_1,\ldots ,X_n)[D_1,\ldots , D_n]$ the al\-ge\-b\-ra of
differential operators. Recall that $A\subset B$ and hence relations  \rf{7}
are
satisfied. Further, each element $f\in B$ can be uniquely represented in the
form
\[
f=\sum_{j_1,\ldots , j_n\ge 0}f_{j_1,\ldots , j_n}D_1^{j_1}\ldots D_n^{j_n}=
\sum_{j\in{\mathbb Z}^n_+}f_jD^j,
\]
where all $f_{j_1,\ldots , j_n}=f_j\in F(X_1,\ldots ,X_n)$ and $F(X_1,\ldots
,X_n)$ is a field of
rational functions over $F$.
Let us replace everywhere in \srf{1} and \srf{2} $A$,  $X^iD^j$,
$\deg f=\deg_{X_1,\ldots ,X_n,D_1,\ldots , D_n}f$, $\dim_FM$, $e_{v,i,j}$,
$f_{v,i,j}\in F$,
$(v,i,j)$, $(i,j)$, $(i',j')$, $(i'',j'')$
by $B$,  $D^j$, $\deg f=\deg_{D_1,\ldots , D_n}f$,  $\dim_{F(X_1,\ldots
,X_n)}M$,
$e_{v,j}$,  $f_{v,j}\in
F(X_1,\ldots ,X_n)$,
$(v,j)$, $j$, $j'$, $j''$ respectively. Thus, we get the definition of the
Janet basis and all
other objects from \srf{1} for the case of the
al\-ge\-b\-ra of differential operators.

We define the homogenization ${^h\!B}$ of $B$ similarly
to ${^h\!A}$, see \srf{3}.
Namely, ${^h\!B}=F(X_1,\ldots ,X_n)[X_0,D_1,\ldots , D_n]$ given by the
relations
\bq{63}
\begin{array}{l}
X_iX_j=X_jX_i,\;D_iD_j=D_jD_i,\quad\mbox{for all}\quad
 i,j, \\
D_iX_i-X_iD_i=X_0,\; 1\le i\le n,\quad X_iD_j=D_jX_i\quad\mbox{for
all}\quad
 i\ne j.\\
\end{array}
\end{equation}
%metka 63 \\
Further, the considerations are similar to the case of the Weyl
al\-ge\-b\-ra $A$ with minor
changes. We leave them to the reader. For example, \trf{2} for the case of
the al\-ge\-b\-ra of differential operators is the same. One need only
to replace everywhere
in its statement $A$, ${^h\!A}$ and  $X_n$ by $B$, ${^h\!B}$ and $D_n$
respectively.
Thus, one can prove \trf{1} for the case when $A$ is
an al\-ge\-b\-ra of differential operators (but now it is $B$).
\trf{1} is proved completely.

One can consider more general algebra of differential operators. Let $\cal F$
be a field with $n$ derivatives $D_1,\dots,D_n$. Then $K_n={\cal F}
[D_1,\dots,D_n]$
is the algebra of differential operators and similarly one can define its
homogenization $^hK_n$ by means of adding the variable $X_0$ satisfying the
relations

$$D_iD_j=D_jD_i,\quad X_0D_i=D_iX_0,\quad D_if-fD_i=f_{D_i}X_0$$

\noindent
for all $i,j$ and any element $f\in \cal F$ where $f_{D_i}\in \cal F$
denotes the
result of
the application of $D_i$ to $f$. Following the proof of \trf{1} one can
deduce
the following statement.

\br{94}
A similar bound to \trf{1} holds for $K_n$.
\end{rems}

\section*{Appendix 1: Degrees of generators
of a graded mo\-du\-le over a polynomial ring and its Hilbert function.}

\noindent We give a short proof of the following result, cf. \cite{Bayer},
\cite{Mora}, \cite{Giusti}, \cite{Dube}.

\bl{14} Let $I\subset {\mathcal A}^l$ be a graded submodule
over the graded polynomial ring ${\mathcal A}=F[X_0,\ldots , X_n]$, and $I$
is
given by a system of
generators $f_1,\ldots ,f_m$ of degrees less than $d$.
Then the Hilbert function $H({\mathcal A}^l/I,m)=\dim_F({\mathcal A}^l/I)_m$
is
stable for $m\ge (dl)^{2^{O(n+1)}}$.
Further, all the coefficients of the Hilbert polynomial
of ${\mathcal A}^l/I$ are bounded
from above by $(dl)^{2^{O(n+1)}}$.
\end{lems}
%metka l14 \\
\noindent{\bf PROOF}\quad Denote $M={\mathcal A}^l/I$.
Let $L\in F[X_0,\ldots , X_n]$ be
a linear form in general position.
Denote by $K$ the kernel of the morphism $M\rightarrow M$ of multiplication
to  $L$.
We have $K=\{z\in {\mathcal A}^l\, :\, Lz=\sum_{1\le i\le
m}f_iz_i,\&\,z_i\in {\mathcal A}\}$.
Hence solving a linear system over ${\mathcal A}$, we get that $K$ has a
system of
generators $g_1,\ldots , g_\mu$ with degrees bounded from above by
$(dl)^{2^{O(n+1)}}$.
Let ${\mathfrak P}$ be an arbitrary associated prime ideal of the mo\-du\-le
$M$ such that ${\mathfrak P}\ne(X_0,\ldots , X_n)$.
Since $L$ is in general position we have  $L\not\in{\mathfrak P}$.
Hence  ${\mathfrak P}$ is not an associated prime ideal of $K$. Therefore,
$K_N=0$ for all sufficiently big $N$. So
$X_i^Ng_j\in I$ for sufficiently big $N$ and all $i,j$.
Hence $g_j=\sum_{1\le i\le m}y_{j,i}f_i$ where $y_{j,i}\in F(X_i)[X_0,\ldots
, X_n]$.
Solving a linear system over the ring $F(X_i)[X_0,\ldots , X_n]$ we get an
estimation for denominators from $F[X_i]$ of all $y_{j,i}$.
Since all $g_j$ and $f_i$ are homogeneous
we can suppose without loss of generality that all the denominators are
$X_i^N$. Thus, we get an upper bound for
$N$.  Namely, $N$ is bounded from
above by $(dl)^{2^{O(n+1)}}$.

Therefore, the sequence
\bq{64}
0\rightarrow M_m\rightarrow M_{m+1}\rightarrow(M/LM)_{m+1}\rightarrow 0
\end{equation}
%metka 64 \\
 is exact
for $m\ge (dl)^{2^{O(n+1)}}$. But $M/LM={\mathcal A}^l/(I+L{\mathcal A}^l)$
is
a
mo\-du\-le over a
polynomial
ring of $F[X_0,\ldots, X_n]/(L)\simeq F[X_0,\ldots, X_{n-1}]$.
Hence by the inductive assumption  the Hilbert function $H({\mathcal
A}^l/(I+L{\mathcal A}^l),m)$ is
stable for
$m\ge (dl)^{2^{O(n)}}$. Therefore, \rf{64} implies that the
Hilbert function $H({\mathcal A}^l/I,m)$ is stable for $m\ge
(dl)^{2^{O(n+1)}}$.

Obviously for $m<(dl)^{2^{O(n+1)}}$ the values $H({\mathcal A}^l/I,m)$ are
bounded from
above
by $(dl)^{2^{O(n+1)}}$.
Hence by the Newton interpolation all the coefficients of the Hilbert
polynomial of ${\mathcal A}^l/I$ are
bounded from above by $(dl)^{2^{O(n+1)}}$.
The lemma is proved.

\medskip\noindent We need also a conversion of \lrf{14}.

\bl{16} Let $I\subset {\mathcal A}^l$ be a graded submodule
over the graded polynomial ring ${\mathcal A}=F[X_0,\ldots , X_n]$.
Assume that the Hilbert function $H({\mathcal A}^l/I,m)=\dim_F({\mathcal
A}^l/I)_m$
is stable for $m\ge D$ and all absolute values  of the coefficients of the
Hilbert polynomial of
the mo\-du\-le ${\mathcal A}^l/I$
are bounded from above by $D$ for some integer $D>1$.
Then $I$ has a system of generators $f_1,\ldots ,f_m$ with degrees
$D^{2^{O(n+1)}}$.
\end{lems}
%metka l16 \\
\noindent{\bf PROOF}\quad Let us choose $f_1,\ldots ,f_m$ to be
the reduced Gr\"obner
basis of $I$ with respect to an admissible linear
order $<$ on the monomials from ${\mathcal A}^l$, cf.
the definitions from \srf{1} and \srf{10}. The degree of a monomial from
${\mathcal A}^l$ is defined similarly to \srf{1} and \srf{10}.
We shall suppose additionally that the considered linear order
is degree compatible,
i.e., for any two monomials $z_1,z_2$ if $\deg z_1<\deg z_2$ then
$z_1<z_2$.
For every $z\in{\mathcal A}$ the greatest monomial
$\mathop{\rm Hdt}\nolimits(z)$ is defined. Further the monomial ideal
$\mathop{\rm Hdt}\nolimits(I)$ is generated by all
$\mathop{\rm Hdt}\nolimits(z)$, $z\in I$. Now $\mathop{\rm
Hdt}\nolimits(f_1),\ldots ,\mathop{\rm Hdt}\nolimits(f_m)$ is a minimal
system
of generators
of $\mathop{\rm Hdt}\nolimits(I)$ and $\deg f_i=\deg\mathop{\rm
Hdt}\nolimits(f_i)$ for every $1\le i\le m$.
The values of Hilbert functions $H({\mathcal A}^l/\mathop{\rm
Hdt}\nolimits(I),m)=H({\mathcal A}^l/I,m)$ coincide for all $m\ge 0$.
Thus, replacing $I$ by $\mathop{\rm Hdt}\nolimits(I)$ we shall assume in
what follows in the proof
that $I$ is a monomial mo\-du\-le.

For every $1\le i\le l$ denote by ${\mathcal A}_i\subset{\mathcal A}^l$
the $i$-th direct summand of
${\mathcal A}^l$. Put $I_i=I\cap{\mathcal A}_i$, $1\le i\le l$. Then
$I\simeq\oplus_{1\le i\le l}I_i$ since $I$ is a monomial mo\-du\-le.
Further, for every $1\le\alpha\le m$ there is $1\le i\le l$ such that
$f_\alpha\in I_i$.
Let us identify ${\mathcal A}_i={\mathcal A}$. Then
$I_i\subset{\mathcal A}$ is a homogeneous monomial ideal. The case
$I_i={\mathcal A}$ is not excluded for some $i$.
For the Hilbert functions we have
\bq{101}
H({\mathcal A}^l/I,m)=\sum_{1\le i\le l}H({\mathcal A}/I_i,m),\quad m\ge 0.
\end{equation}
%metka 101 \\
If $({\mathcal A}/I_i)_D=0$ for some $i$
then $({\mathcal A}/I_i)_m=0$ for every $m\ge
D$.
In this case the ideal $I_i$ is generated by $\sum_{0\le m\le D}(I_i)_m$.
Hence in \rf{101} for the values $m\ge D$
one can omit this index $i$ in the sum from the right part.
Therefore, in this case the
proof is reduced to a smaller $l$.
So we shall assume without loss of generality that $({\mathcal A}/I_i)_D\ne
0$, $1\le i\le l$.

Further, we use the exact description of the Hilbert function of a
homogeneous ideal, see \cite{Dube}
Section~7.
Namely there are the unique integers $b_{i,0}\ge b_{i,1}\ge\ldots\ge
b_{i,n+2}=0$ such
that
\bq{100}
H({\mathcal A}/I_i,m)={m+n+1 \choose n+1}-1-
\sum_{1\le j\le n+1}{m-b_{i,j}+j-1 \choose j}
\end{equation}
%metka 100 \\
for all sufficiently big $m$ and
\bq{103}
b_{i,0}=\min\{d\, :\, d\ge
b_{i,1}\,\&\,\forall\,m>d\quad\mbox{\rf{100}}\quad\mbox{holds}\, \}.
\end{equation}
%metka 103 \\
This description (without constants $b_{i,0}$)
is originated from the classical paper \cite{Mac}.
The integers $b_{i,0},\ldots,
b_{i,n+2}$ are called the Macaulay constants of the ideal $I_i$.
Besides that,
\bq{104}
h(i,m)=H({\mathcal A}/I_i,m)-{m+n+1\choose n+1}+1+
\sum_{1\le j\le n+1}{m-b_{i,j}+j-1\choose j}\ge 0
\end{equation}
%metka 104 \\
for every $m\ge b_{i,1}$, see \cite{Dube} Section~7.
By Lemma~7.2 \cite{Dube}
for all $1\le\alpha\le m$ if $f_\alpha\in I_i$ then
$\deg f_\alpha\le b_{i,0}$.
Hence it is sufficient to prove that all $b_{i,0}$, $1\le i\le l$,
are bounded from above by $D^{2^{O(n+1)}}$.

By  \rf{101} and  \rf{100}  the coefficient at $m^{n-j}$, $0\le j\le n$,
of the Hilbert polynomial of ${\mathcal A}^l/I$ is
\bq{102a}
\frac{\mu_j}{(n+1-j)!}\sum_{1\le i\le l}b_{i,n+1-j}+
\sum_{0\le v\le j-1}\sum_{1\le i\le l}
\frac{1}{(n+1-v)!}\mu_{j,v}(b_{i,n+1-v}),
\end{equation}
%metka 102a \\
where $0\ne\mu_j$ is an integer and $\mu_{j,v}\in{\mathbb Z}[Z]$, $0\le v\le
j-1$,
is a polynomial with integer coefficients with $\deg\mu_{j,v}=j-v+1$.
Moreover,
$|\mu_j|$ and
absolute values of all the coefficients
of all the polynomials $\mu_{j,v}$ are
bounded from above by, say,
$2^{O(n^2)}$.
Denote $b_j=\sum_{1\le i\le l}b_{i,j}$, $0\le j\le n+2$.
By the condition of the lemma all the coefficients
of the Hilbert polynomial of ${\mathcal A}^l/I$ are bounded from above by
$D$.
Hence from \rf{102a} one can recursively estimate $b_{n+1},b_n,\ldots , b_1$.
Namely, $b_{n+1-j}=(2^{n^2}lD)^{2^{O(j+1)}}$, $0\le j\le n$.
Hence $b_1=(lD)^{2^{O(n+1)}}$.  Notice that $b_{i,1}\le\max_{1\le i\le
l}{b_{i,1}}\le b_1$
for every $1\le i\le m$.

Now let $m\ge\max_{1\le i\le l}{b_{i,1}}$. By \rf{104} if
$h(i,m)\ne 0$ for some $1\le i\le l$ then $m<D$, i.e., $m$
is less than the bound  $D$
for the stabilization of the Hilbert function of ${\mathcal A}^l/I$. Thus,
$b_{i,0}\le\max\{b_{i,1},D\}$ by \rf{103}.
Hence $b_{i,0}$ is bounded from above by $(lD)^{2^{O(n+1)}}$.

We have $({\mathcal A}/I_i)_D\ne 0$ for every $1\le i\le l$.
This implies $H({\mathcal A}^l/I,D)\ge l$.
Denote by $c_j$ the $j$-th coefficient of the Hilbert polynomial of the
mo\-du\-le  ${\mathcal A}^l/I$.
Now $|c_j|D^j\ge l/(n+1)$ for at least one $j$.
Hence $D^{n+1}(n+1)\ge l$ by the condition of the lemma.
This implies that $l^{2^{O(n+1)}}$ is bounded from above by $D^{2^{O(n+1)}}$.
Therefore, $b_{i,0}$ is bounded from above by $D^{2^{O(n+1)}}$.
The lemma is proved.

\section*{Appendix 2: Bound on the Gr\"obner basis of a monomial module via
the coefficients of its
Hilbert polynomial}

Denote by $C_l=\ZZ_+^n\cup \cdots \cup \ZZ_+^n$ the disjoint union of $l$
copies of the
semigrid $\ZZ_+^n=\{(i_1,\dots,i_n):i_j\geq 0, 1\leq j\leq n\}$.
A subset of $C_l$ which intersects each disjoint copy
of $\ZZ_+^n$ by a semigroup closed with respect to
addition of elements from $\ZZ_+^n$ is called an ideal of $C_l$.
Any ideal $I$ in
$C_l$ has a unique finite Gr\"obner basis $V=V_I$, denote $T=C_l \setminus
I$.
Clearly, $I$
corresponds to a monomial submodule in the free module $(F[X_1,\dots,
X_n])^l$. The
degree of an element $u=(k;i_1,\dots,i_n)\in C_l, 1\leq k\leq l$ is defined
as
$|u|=i_1+\cdots+i_n$. The degree of a subset in $C_l$ is defined as the
maximum of the
degrees of its elements. The Hilbert function $H_T(z)$ equals to the
number of vectors
$u\in T$ such that $|u|\leq z$. Then $H_T(z)=\sum_{0\leq s\leq m} c_sz^s,
\quad
z\geq z_0$
for suitable $z_0$, integers $c_0,\dots,c_m$ where the degree $m\leq n$.
Denote
$c=\max_{0\leq s\leq m} |c_s|s!+1$.

\bprp{monomial}
(cf. \cite{Giusti}, \cite{Mora}, \cite{Dube}). The degree of $V$ does not
exceed
$(cn)^{2^{O(m)}}$.
\end{props}

\noindent{\bf PROOF}\quad An {\it $s$-cone} we call a subset of a $k$-th
copy of
$\ZZ_+^n$ in $C_l$ for
a certain $1\leq k\leq l$ of the form

\begin{equation}\label{90}
P=\{X_{j_1}=i_1,\dots,X_{j_{n-s}}=i_{n-s}\}
\end{equation}

\noindent
for suitable $1\leq j_1,\dots, j_{n-s}\leq n$. The degree of \rf{90} we
define as
$|P|=i_1+\cdots+i_{n-s}$ (note that this definition is
different from the one in \cite{Dube}). By a {\it predessesor} of \rf{90} we
mean each
$s$-cone in the same $k$-th copy of $\ZZ_+^n$ of the
type

\begin{equation}\label{91}
\{X_{j_1}=i_1,\dots,X_{j_{p-1}}=i_{p-1},X_{j_p}=i_p-1,X_{j_{p+1}}=i_{p+1},\dots,
X_{j_{n-s}}=i_{n-s}\}
\end{equation}

\noindent
for some $1\leq p\leq n-s$, provided that $i_p\geq 1$. Fix an arbitrary
linear order
on $s$-cones compatible with the relation of predessesors.

By inverse recursion on $s$ we fill gradually $T$ (as a union) by $s$-cones.
For the base
we start with $s=m$. Assume that a current union $T_0\subset T$ of $m$-cones
is already
constructed (at the very beginning we put $T_0=\emptyset$) and an $m$-cone
of the form
 \rf{90} with $s=m$ is the least one (with respect to the fixed linear order
on
$m$-cones) which is contained in $T$ not being a subset of $T_0$. Observe
that each
predessesor of this $m$-cone was added to $T_0$ at earlier steps of its
construction.
Since the total number of $m$-cones added to $T_0$ does not exceed $c_mm!<
c$ we
deduce that the degree of every such $m$-cone is less than $c_mm!$ (taking
into account
that the very first $m$-cone added to $T_0$ has the degree $0$).

For the recursive step assume that the current $T_0$ is a union of
all possible $m$-cones, $(m-1)$-cones,...,$(s+1)$-cones and
perhaps, some $s$-cones. This can be expressed as
$\deg(H_T-H_{T_0})\leq s$. Again as in the base take the least
$s$-cone of the form \rf{90} which is contained in $T$ not being a
subset of $T_0$. Observe that each predessesor of the type \rf{91}
of this $s$-cone is contained in an appropriate $r$-cone $Q$,
$r\ge s$, such that $Q$ was added to $T_0$ at earlier steps of its
constructing and $Q\subset\{X_{j_p}=i_p-1\}$. Hence

\begin{equation}\label{92}
%\{X_{j_{l(1)}}=i'_{l(1)},\dots, X_{j_{l(n-r)}}=i'_{l(n-r)}\}
|Q|\ge i_p-1.
\end{equation}

%\noindent
%with the extra property that $j_p$ is among $j_{l(1)},\dots,j_{l(n-r)}$ for
%a certain
%$r\geq s$, and \rf{92} was added to $T_0$ at earlier steps of its
%constructing. Thus,
%$i'_{l(q)}=i_{l(q)}$ for all $1\leq q\leq n-r$ except for $q_0$ such that
%$l(q_0)=p$,
%and in the latter case we have $i'_p=i_p-1$. In the sequel, this extra
%property would
%allow one to estimate $i_p$ (and similarly the degree of \rf{90} considering
%other
%predessesors of \rf{90}) by inverse induction
%on $s$.
The described construction terminates when $T_0=T$.
Denote by $t_s$ the number of $s$-cones added to $T_0$ and by $k_s$ the
maximum of their
degrees. We have seen already that $t_m,k_m< c$.

%Exploiting repeatedly the
%inclusion-exclusion formula one can verify
%by inverse induction on $s$ that
%$t_s,k_s\leq c^{n^{O(m-s)}}$, taking into account
%the stated property of the
%cones of the
%form \rf{91}.

Now by inverse induction on $s$ we prove that $t_s,k_s\leq
(cn)^{2^{O(m-s)}}$. To this end
we introduce a relevant semilattice on cones. Let ${\cal
C}=\{C_{\alpha,\beta}\}_
{\alpha,\beta}, \quad 0\leq \beta \leq \gamma_{\alpha}$ be a family of cones
of the form
\rf{90} where $\dim C_{\alpha,\beta}=\alpha$. By an $\alpha$-piece we call an
$\alpha$-cone being the intersection of a few cones from ${\cal C}$. All the
pieces
constitute a semilattice $\cal L$ with respect to the intersection and with
maximal
elements from ${\cal C}$. We treat $\cal L$ also as a partially ordered set
with
respect to the inclusion relation. Clearly, the depth of $\cal L$ is less
than $n$.
Our nearest purpose is to bound from above  the size of $\cal L$. For the
sake of
simplifying the bound we assume (and this will suffice for our goal in the
sequel) that
$\gamma_{\alpha}\leq (cn)^{2^{O(m-\alpha)}}$ for $s\leq \alpha \leq m$ and
$\gamma_{\alpha}
=0$ when $\alpha <s$, although one could write a bound in general in the
same way. Besides that we assume that the constant in $O(\ldots)$ is
sufficiently big.
In what follows all the constants in $O(\ldots)$ coincide.

\bl{size}
Under the assumption on the numbers
$\gamma_{\alpha}\leq (cn)^{2^{O(m-\alpha)}}, \quad s\leq \alpha \leq m$ of
maximal elements of all
dimensions
from ${\cal C}$,
the number of $\alpha$-pieces in $\cal L$ does not exceed
$(cn)^{2^{O(m-\alpha)}+1}$ for
$s\leq \alpha \leq m$ or $(cn)^{2^{O(m-s)}(s-\alpha+1)+1}$ when $\alpha <s$.
\end{lems}

\noindent{\bf PROOF}\quad For each $\alpha$-piece choose its arbitrary
irredundant representation
as the intersection of the cones from ${\cal C}$. Let $\delta$ be the minimal
dimension among these cones. Then this intersection contains at most
$\delta-\alpha+1$
cones. Therefore, the number of possible $\alpha$-pieces does not exceed
$$\sum_{max\{\alpha,s\}\leq \delta \leq m}
(cn)^{2^{O(m-\delta)}(\delta-\alpha+1)},$$
\noindent
that proves the lemma.

Now we come back to estimating $t_s,k_s$ by inverse induction on $s$. Let in
the
described above construction the current $T_0$ is the union of all added
$m$-cones,
$(m-1)$-cones,...,$s$-cones. Denote this family of cones by ${\cal C}$ and
consider
the corresponding semilattice $\cal L$ (see above). Our next purpose is to
represent
$T_0$ as a $\ZZ$-linear combination of the pieces from $\cal L$ by means of
a kind of
the inclusion-exclusion formula. We assign the coefficients of this
combination by
recursion in $\cal L$. As a base we assign 1 to each maximal piece, so to
the elements
of ${\cal C}$. As a recursive step, if for a certain piece $P\in \cal L$ the
coefficients are already assigned to all the pieces greater than $P$, we
assign to
$P$ the coefficient $\epsilon_P$ in such a way that the sum of the assigned
coefficients to $P$ and to all the greater pieces equals to 1. Therefore, we
get
$$T_0=\sum_{P\in \cal L} \epsilon_P P$$
\noindent
where the sum is understood in the sense of multisets. Hence

\begin{equation}\label{102}
H_{T_0}(z)=\sum_{P\in \cal L} \epsilon_P {z-|P|+\dim P \choose \dim P}
\end{equation}

\noindent
for large enough $z$. We recall that $\deg(H_T-H_{T_0})\leq s-1$.

Now we majorate the  coefficients $|\epsilon_P|$ by induction in the
semilattice
$\cal L$. The inductive hypothesis on $t_{\alpha}\leq (cn)^{2^{O(m-\alpha)}},
s\leq \alpha \leq m$ and \lrf{size} imply that
$$\sum_{dim P=\lambda} |\epsilon_P|\leq (cn)^{2^{O(m-\lambda)}}, \quad
s-1\leq \lambda \leq m.$$
\noindent
by inverse induction on $\lambda$ following the assigning $\epsilon_P$.
In fact, one could majorate in a similar way also $\sum_{\dim P=\lambda}
|\epsilon_P|$
when $\lambda < s-1$, but we don't need it. The inductive hypothesis on
$k_{\alpha}\leq (cn)^{2^{O(m-\alpha)}}, \quad s\leq \alpha \leq m$ and
\rf{102} entail that
the coefficient of $H_{T_0}(z)$ at the power $z^{\alpha}$ does not exceed
$(cn)^{2^{O(m-\alpha)}}, \quad s-1\leq \alpha \leq m$ (actually, due to the
inequality
$\deg(H_T-H_{T_0})\leq s-1$ the coefficients at the powers $z^{\alpha}$ for
$s\leq \alpha \leq m$ are less than $c$). In particular, the coefficient at
the power
$z^{s-1}$ does not exceed $(cn)^{2^{O(m-s+1)}}$. Denote $H_T-H_{T_0}=\eta
z^{s-1}+\cdots$.
By constructing $T_0$ we add to it $t_{s-1}=\eta (s-1)!$ of $(s-1)$-cones,
which
justifies the inductive step for $t_{s-1}\leq (cn)^{2^{O(m-s+1)}}$.

To conduct the inductive step for $k_{s-1}\leq (cn)^{2^{O(m-s+1)}}$ we
observe that for each $(s-1)$-cone $P$ added to $T_0$ either every its
predessesor
%of the form \rf{92}
is contained in
a cone of dimension at least $s$, or some its predessesor is an $(s-1)$-cone
as well. In
the former case $|P|\leq (\max_{s\leq \alpha \leq m}k_{\alpha}+1)(n-s+1)$
(due to \rf{92}), while
%the stated property just after \rf{92}), while
in the latter case $|P|$ is greater by 1 than the degree of this predessesor,
hence $k_{s-1}\leq (\max_{s\leq \alpha \leq m}k_{\alpha}+1)(n-s+1)+t_{s-1}$.
Finally,
exploit  the inductive hypothesis for $k_m,\dots,k_s$, and the just obtained
inequality on $t_{s-1}$.

To complete the proof of the proposition it suffices to notice that for any
vector from the
basis $V$ treated  as an $0$-cone, each its predessesor of the type \rf{91}
for $s=0$ is
contained in an appropriate $r$-cone,
whence the degree of $V$ does not exceed
$(\max_{0\leq \alpha \leq m}k_{\alpha}+1)n$ again due to \rf{92} (cf. above).

\medskip{\bf Acknowledgement.} The authors are grateful to the Max-Planck
Institut f\"ur
Mathematik, Bonn for its hospitality during the stay where the paper was
written.

%\newpage

\end{document}